\documentclass[11pt]{article}
\usepackage{graphicx}
\usepackage{epstopdf}
\usepackage{amsmath, amssymb}
\usepackage{latexsym, color}
\usepackage{mathtools}
\usepackage{xcolor}
\usepackage{graphicx}

\def\la{\big\langle}
\def\ra{\big\rangle}

\def\ds{\displaystyle}
\def\forall{\hbox{for all}~}
\def\L{{\bf L}}

\def\bfv{{\bf v}}
\def\bfw{{\bf w}}

\def\bfn{{\bf n}}

\def\ve{\varepsilon}

\def\E{{\cal E}}
\def\A{{\cal A}}

\def\H{{\cal H}}

\def\avint{-\!\!\!\!\!\!\int}

\def\R{I\!\!R}

\def\vp{\varphi}

\def\vs{\vskip 2em}

\def\v{\vskip 1em}
\def\O{{\cal O}}

\def\C{{\cal C}}
\def\J{{\cal J}}

\def\H{{\cal H}}
\def\rest{\llcorner}
\def\bega{\begin{array}}
\def\enda{\end{array}}
\def\begi{\begin{itemize}}
\def\endi{\end{itemize}}
\def\ov{\overline}
\def\wto{\rightharpoonup}
\def\Tilde{\widetilde}
\def\Hat{\widehat}

\def\meas{\hbox{meas}}
\def\bel{\begin{equation}\label}
\def\eeq{\end{equation}}
\def\sqr#1#2{\vbox{\hrule height .#2pt
\hbox{\vrule width .#2pt height #1pt \kern #1pt
\vrule width .#2pt}\hrule height .#2pt }}
\def\square{\sqr74}
\def\endproof{\hphantom{MM}\hfill\llap{$\square$}\goodbreak}
\definecolor{cadmiumgreen}{rgb}{0.0, 0.42, 0.24}

\DeclareMathOperator*{\argmin}{arg\,min}

\newtheorem{theorem}{Theorem}[section]

\newtheorem{lemma}{Lemma}[section]

\newtheorem{remark}{Remark}[section]

\newtheorem{example}{Example}[section]
\newtheorem{conjecture}{Conjecture}[section]

\begin{document}
\title{\bf Optimal Control of Moving Sets}\vs
\author{Alberto Bressan, Maria Teresa Chiri, and Najmeh Salehi\\
\,
\\
Department of Mathematics, Penn State University \\
University Park, Pa.~16802, USA.\\
\,
\\
e-mails: axb62@psu.edu, mxc6028@psu.edu, nfs5419@psu.edu.
}
\maketitle

\begin{abstract} Motivated by the  control of invasive biological populations,
we consider a class of optimization problems for moving sets $t\mapsto \Omega(t)\subset\R^2$.
Given an initial set $\Omega_0$, the goal is to minimize the area of the contaminated set 
$\Omega(t)$ over time, 
plus a cost related to the control effort.  
Here the control function is the inward normal speed along the boundary $\partial \Omega(t)$. We prove the existence of optimal solutions, within a class of sets with finite perimeter.   Necessary conditions for optimality
are then derived, in the form of a Pontryagin maximum principle.  Additional optimality conditions show that the
sets $\Omega(t)$ cannot have certain types of outward or inward corners. 
Finally, some explicit solutions are presented.
\end{abstract}

\section{Introduction}
\label{s:1}
\setcounter{equation}{0}
The original motivation for our study comes from control
problems for reaction-diffusion equations, of the form
\bel{CRD}
u_t~=~f(u) + \Delta u - g(u, \alpha).\eeq
Here we think of $u=u(t,x)$ as the density of an invasive population, which
grows at rate $f(u)$, diffuses in space, and can by partly removed by
implementing a control $\alpha=\alpha(t,x)$.  Given an initial 
density 
$$u(0,x)~=~\bar u(x),$$
we seek
a control $\alpha$ which minimizes:
$$\hbox{[total size of the population over time] + [cost of implementing the control].}$$
For example, $u$ may describe the density of mosquitoes, the control
$\alpha$ is the amount of pesticides which are used, and $g(u,\alpha)$ is the rate at which 
mosquitoes are eliminated by this action.
Assuming that the source term has two 
equilibrium states at $u=0$ and $u=1$, i.e.
$$f(0)~=~f(1)~=~0,$$
a simplified model was derived in \cite{BCS},
formulated in terms of the evolution of a set.  Indeed, in a common situation one can identify
a ``contaminated set" 
$$\Omega(t)~=~\bigl\{x\in \R^d\,;~u(t,x)\approx 1\bigr\}$$
where the population is large,
while $u(t,x)\approx 0$ over most of the complementary set $\R^d\setminus \Omega(t)$.
By implementing different controls $\alpha$ in (\ref{CRD}), one can shrink the 
contaminated set $\Omega(t)$ at different rates. 

As shown by the analysis in \cite{BCS}, the effort $E(\beta)$ required for pushing the boundary 
$\partial \Omega(t)$ inward, with speed $\beta$ in the normal direction,
can be determined in terms of a minimization problem.   
Namely, we define $E(\beta)$ to be the minimum cost of a control $\alpha$ which yields 
a traveling wave profile for (\ref{CRD}) with speed $\beta$. 
More precisely,
$$E(\beta)~=~\min \|\alpha\|_{\L^1}\,,$$ where the minimization is taken over all controls 
$\alpha:\R\mapsto\R_+$ such that there exists a profile $U:\R\mapsto [0,1]$
with
$$U'' + \beta  U' + f(U) - g(U,\alpha)~=~0, \qquad \qquad U(-\infty) = 0,\qquad U(+\infty) =1.$$
By taking a sharp interface limit (see again \cite{BCS}) for details),
this leads to a class of optimization problems for the evolution of a set $t\mapsto \Omega(t)$
of finite perimeter.   Here we regard the inward normal velocity  $\beta=\beta(t,x)$, 
assigned at every point $x\in \partial \Omega(t)$, as our new control function.
The main goal of the present paper is to study these optimization problems for a moving set.
For notational simplicity, we carry out the analysis in the planar case $\Omega(t)\subset \R^2$,
which is most relevant for applications.  
We expect that our results can be extended to higher space dimensions, with similar proofs.

Three problems will be considered.
\begi
\item[{\bf (NCP)}] {\bf Null Controllability Problem.}  
{\it   Let an initial set $\Omega_0\subset\R^2$,  a convex cost 
function $E:\R\mapsto\R_+$, and a constant $M>0$
be given. 
Find a set-valued function $t\mapsto \Omega(t)$  such that, for some $T>0$,
\bel{nco}\Omega(0)~=~\Omega_0,\qquad\qquad 
\Omega(T)~=~\emptyset,\eeq
\bel{eb}\E(t)~\doteq~
\int_{\partial\Omega(t)}E\bigl(\beta(t,x)\bigr)\, d\sigma~\leq~ M\quad \qquad \forall t\in [0,T],\eeq
where  $\beta$ denotes the velocity of a boundary point in the inward normal direction, 
and the integral is taken w.r.t.~the 1-dimensional Hausdorff measure 
along the boundary of $\Omega(t)$.}
\v
\item[{\bf (MTP)}]  {\bf Minimum Time Problem}  {\it  Among all strategies that satisfy
(\ref{nco})-(\ref{eb}), find one which minimizes the time $T$.}
\v
\item[{\bf (OP)}] {\bf Optimization Problem.}  {\it Let an initial set
$\Omega_0\subset\R^2$ 
and cost functions  $E: \R\mapsto \R_+$, $\phi:\R_+\mapsto\R_+\cup\{+\infty\}$ be given.
Find a set-valued function $t\mapsto \Omega(t)$, with $\Omega(0)=\Omega_0$,
which minimizes
\bel{cost1}
J~=~\int_0^T\phi\bigl(\E(t)\bigr) dt + 
c_1\int_0^T m_2\bigl(\Omega(t)\bigr)\, dt + c_2\,m_2\bigl(\Omega(T)\bigr).\eeq
Here $\E(t)$ is the total effort at time $t$, defined as in (\ref{eb}), while $m_2$ denotes 
the 2-dimensional Lebesgue measure.
}
\endi
%

\begin{remark}{\rm In (\ref{cost1}), we are thinking of $\Omega(t)\subset\R^2$ as a contaminated region.  We seek to minimize the area of this region,
plus a cost related to the control effort. 
}\end{remark}

We now introduce assumptions on the cost functions
 $E$ and $\phi$, that will be used in the sequel.
\begi
\item[{\bf (A1)}]  {\it The function $E:\R\mapsto\R_+$ is continuous and convex.
There exist constants $\beta_0<0$ and $a>0$ such that
\bel{Eprop}
\left\{ \bega{rll}
E(\beta)&\geq ~a(\beta-\beta_0)\qquad &\hbox{if}~~\beta\geq \beta_0,\\[2mm]
E(\beta)&=~0\qquad &\hbox{if}~~\beta\leq \beta_0.\enda\right.
\eeq
In addition,  $E$ is twice continuously differentiable for $\beta>\beta_0$ and satisfies
\bel{Eass} E(\beta) - \beta\, E'(\beta)~\geq~0\qquad\qquad\forall \beta>0.\eeq
}
\endi
\begi
\item[{\bf (A2)}]  {\it The function $\phi:\R_+\mapsto\R_+
\cup\{+\infty\}$  is lower semicontinuous, nondecreasing, and convex.
Moreover, for some constants $C_1, C_2>0$ and $p>1$ one has
\bel{phia}
\phi(0)=0, \qquad \quad \phi(s)~\geq~C_1 s^p - C_2\qquad\forall ~s\geq 0\,.\eeq
}
\endi
\begin{remark}\label{r:12} {\rm 
The particular form (\ref{Eprop}) of the effort function $E$ models the fact that, if no 
control is applied, the contaminated region $\Omega(t)$ expands in all directions
with speed $|\beta_0|$. If some control is active, this expansion rate can be reduced, or even 
reversed (so that the contaminated region actually shrinks).   
For practical purposes, the values of the function $E$ 
for inward normal speed $\beta<\beta_0<0$ are irrelevant, 
because in an optimization problem it is never convenient to let the set $\Omega(t)$ 
expand with speed larger than $|\beta_0|$. Indeed, this will only increase the 
overall cost in (\ref{cost1}).   As it will be shown in Section~\ref{s:2}, thanks to
the choice $E(\beta)=0$ in (\ref{Eprop})  one achieves the convexity of the 
cost function $L$ in (\ref{Lext}).
}
\end{remark}

\begin{remark}\label{r:13} {\rm
The assumptions {\bf (A1)} imply that the effort function 
$E$ has sublinear growth. Indeed, when $\beta> 1$,  from (\ref{Eass}) by convexity
it follows
$$E(\beta)- \beta \,{E(\beta)-E(1)\over\beta-1}~\geq~E(\beta)- \beta \,E'(\beta)~\geq ~0.$$
This implies
\bel{Esg}
E(\beta)~\leq~\beta\, E(1)\qquad\qquad \forall \beta\geq 1.\eeq
}
\end{remark}

\begin{figure}[ht]
\centerline{\hbox{\includegraphics[width=16cm]{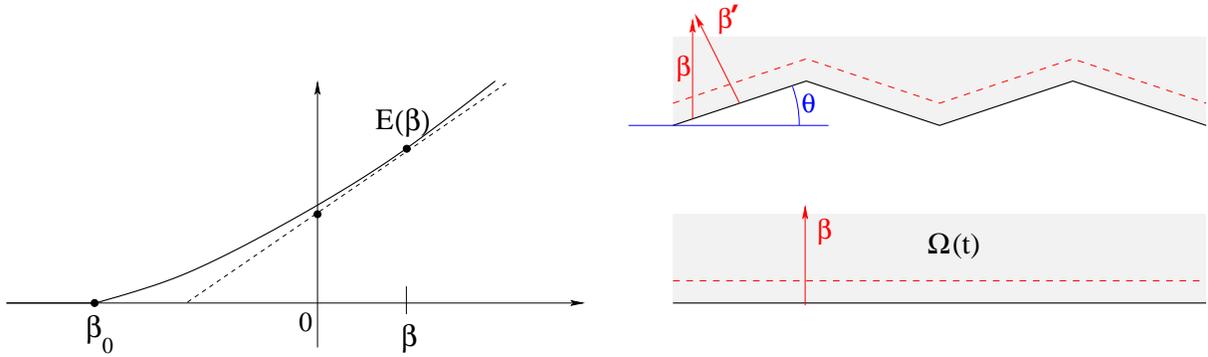}}}
\caption{\small Left: an effort function $E$ that satisfies the assumption (\ref{Eass}).
Right: a geometric explanation of the necessity of this assumption.
}
\label{f:sc54}
\end{figure}

\begin{remark}\label{r:14}
{\rm  The geometric meaning of the assumption (\ref{Eass}) is shown in Fig.~\ref{f:sc54}, left.
Namely, the tangent line at $\beta$ has a positive intersection with the vertical axis.
This assumption  plays a crucial role for the lower semicontinuity of the 
cost functional in (\ref{cost1}). 
Indeed, consider the 
situation in Figure~\ref{f:sc54}, right. Here the contaminated region $\Omega(t)$ shrinks with 
inward normal speed $\beta>0$.  We can slightly perturb its boundary, so that it 
oscillates with angle $\pm \theta$.
In the first case the cost of the control, for a portion of the boundary  parameterized by $s\in [s_1, s_2]$,  is computed by
$$\int_{s_1}^{s_2} E(\beta)\, ds,$$
where $ds$ denotes arclength along the boundary.
In the perturbed case, the boundary length increases, while the normal speed decreases.
This leads to
$$ds'~ =~ {ds\over \cos\theta}\,,\qquad\qquad \beta' ~=~\beta\, \cos\theta.$$
The new total effort is
$$\int_{s_1}^{s_2} {E(\beta \cos\theta)\over \cos\theta}\, ds\,.$$
Differentiating w.r.t.~$\theta$ we find
$${d\over d\theta}  {E(\beta \cos\theta)\over \cos\theta}\Bigg|_{\theta=0}~=~0,\qquad\qquad
{d^2\over d\theta^2}  {E(\beta \cos\theta)\over \cos\theta}\Bigg|_{\theta=0}~=~E(\beta) -
\beta\, E'(\beta).$$
Therefore, if we want the wiggly profile to yield a larger or equal cost, 
the inequality (\ref{Eass}) must be satisfied.
}
\end{remark}

\begin{example}\label{e:11}
{\rm 
Consider the
basic case where
\bel{EP}\left\{ \bega{rll}  E(\beta)&=~\beta+1\qquad &\hbox{if}~~\beta\geq -1,\cr
E(\beta)&=~0\qquad &\hbox{if}~~\beta\leq -1,\enda\right.
\qquad\qquad 
 \phi(s)~=~\left\{ \bega{cl} 0\quad &\hbox{if}\quad s\leq M,\cr
+\infty\quad &\hbox{if}\quad s> M.\enda\right.\eeq
In this case, with zero control effort, we obtain a set $\Omega(t)$ 
which expands in all directions
with unit speed.
On the other hand,  the bound (\ref{eb}) on  the instantaneous control effort 
allows us reduce the area of 
$\Omega(t)$ at rate $M$ per unit time.
Calling 
$$A(t)~=~m_2(\Omega(t)), \qquad P(t)~=~m_1(\partial\Omega(t))$$
respectively the area (i.e., the 2-dimensional Lebesgue measure) of $\Omega(t)$
 and the perimeter of 
$\Omega(t)$, (i.e., the 1-dimensional Hausdorff measure
of the boundary $\partial \Omega(t)$),  
we thus have
\bel{AS}{d\over dt} A(t)~=~P(t) -M.\eeq
The Null Controllability Problem {\bf (NCP)} can thus be solved if we can reduce the perimeter $P(t)$
to a value smaller than $M$.
}
\end{example}

The remainder of the paper is organized as follows. 
Section~\ref{s:2} contains preliminary material.   In particular, we reformulate
the optimization problems in terms of sets with finite perimeter, and prove
the one-sided H\"older estimate (\ref{1sh}) 
on the area of the sets $\Omega(t)$, valid whenever the 
total cost of the control is bounded.
In Section~\ref{s:3} we give a simple condition that ensures the solvability of the 
Null Controllability Problem.
The existence of optimal solutions to {\bf (OP}) is proved in Section~\ref{s:4},
in the somewhat relaxed formulation (\ref{FU}).   
The following Sections~\ref{s:5} to~\ref{s:7} establish various necessary 
conditions for optimality.
Finally, in Section~\ref{s:8} we discuss the geometric meaning of the necessary conditions,
and give an example of a set-valued map
$t\mapsto \Omega(t)$ that satisfies all these necessary conditions.   

Based on these necessary conditions,
we expect that at each time $t\in [0,T]$ the optimal control should concentrate all the effort 
along the portion of the boundary $\partial \Omega(t)$
with maximum curvature.
At the present time, however,  this remains on open question (see Conjecture~\ref{c:81}), for two reasons:
\begi
\item[(i)] The existence of optimal solutions is here proved within a class of sets with finite perimeter,
while our necessary conditions require piecewise $\C^2$ regularity.
\item[(ii)] The set-valued functions $t\mapsto \Omega(t)$ 
which satisfy our optimality conditions may not be unique.
\endi

In earlier literature, several different models related to the control of a moving set have been analyzed in 
\cite{Bblock, Breview, BMN, BZ, CPo, CLP}.    
Eradication problems for invasive biological 
species were studied in \cite{ACD, ACM}. 

\section{Preliminaries}
\label{s:2}
\setcounter{equation}{0}
The instantaneous effort functional $\E(t)$ in (\ref{eb}) is naturally defined for 
moving sets with $\C^1$ boundary.    Toward the analysis of optimization problems, 
it will be convenient to extend this definition to more general sets with finite perimeter
\cite{AFP, M}.
We shall thus work within the family of admissible sets
\bel{Adef} \A~\doteq~\Big\{  \Omega\subset \,]0,T[\,\times\R^2\,;~~\Omega ~\hbox{is bounded
and has bounded perimeter}\Big\}.\eeq  
Calling ${\bf 1}_\Omega$ the characteristic function of $\Omega$, 
this implies that  ${\bf 1}_\Omega\in BV$.
In other words, the distributional gradient $\mu_\Omega\doteq D\, {\bf 1}_\Omega$
is a finite $\R^3$-valued Radon measure:
\bel{div}\int_\Omega \hbox{div}\, \vp\, dx~=~-\int \vp\cdot d\mu_\Omega\qquad\qquad \forall \vp\in\C^1_c
\bigl(]0,T[\,\times\R^2\,;~\R^3\bigr).\eeq
Given a set $\Omega\in \A$, we consider the multifunction
\bel{slice} t~\mapsto~\Omega(t)~\doteq~\bigl\{ x\in\R^2\,;~(t,x)\in \Omega\bigr\}.\eeq
By possibly modifying ${\bf 1}_\Omega$ on a set of 3-dimensional measure zero,
the map $t\mapsto {\bf 1}_{\Omega(t)}$
has bounded variation from $\,]0,T[\,$ into $\L^1(\R^2)$.  In particular, for every $0<t<T$,
the one-sided limits
\bel{ompm}{\bf 1}_{\Omega(t+)}~\doteq~\lim_{t\to t+} {\bf 1}_{\Omega(t)}\,,
\qquad\qquad {\bf 1}_{\Omega(t-)}~\doteq~\lim_{t\to t-} {\bf 1}_{\Omega(t)}\,,\eeq
are well defined in $\L^1(\R^2)$.   This uniquely defines the sets $\Omega(t+)$,
$\Omega(t-)$, up to a set of 2-dimensional Lebesgue measure zero.
Throughout the following, we define
the sets $\Omega(0)$ and $\Omega(T)$ in terms of
\bel{om0T}{\bf 1}_{\Omega(0)}~\doteq~\lim_{t\to 0+} {\bf 1}_{\Omega(t)}\,,
\qquad\qquad {\bf 1}_{\Omega(T)}~\doteq~\lim_{t\to T-} {\bf 1}_{\Omega(t)}\,.\eeq
We  write $\H^m$ for the $m$-dimensional Hausdorff measure,
while $\mu\,\rest V$ denotes the restriction of a measure $\mu$ to the set $V$.
By $B(y,r)$ we denote the open ball centered at $y$ with radius $r$, while
$S^2$ is the sphere of unit vectors in $\R^3$.
%
%
%
%

For every set of finite perimeter 
$\Omega\in \A$, its reduced boundary  $\partial^*\Omega$ is defined to be the set of points
$y=(t,x)\in\,]0,T[\,\times\R^2$ such that 
\bel{nuo}
\nu_\Omega(y)~\doteq~\lim_{r\downarrow 0}~ {\mu_\Omega\bigl(B(y,r)\bigr)\over
|\mu_\Omega|\bigl(B(y,r)\bigr)}\eeq
exists in $\R^3$ and satisfies $|\nu_\Omega(y)|=1$. 
The function $\nu_\Omega: \partial^*\Omega\mapsto S^2$ is called the 
{\it generalized inner normal}
to $\Omega$.    A fundamental theorem of De Giorgi \cite{AFP, M} implies that
$\partial^*\Omega$ is countably 2-rectifiable and $|D {\bf 1}_\Omega| = \H^2\rest \partial^*\Omega$.

In order to introduce a cost associated with each set $\Omega\in\A$, 
we observe that, in the smooth case,  the (inward) normal velocity of the set 
$\Omega(t)$  at the point 
$(t,x)\in \partial^*\Omega$ is computed by
\bel{bdef}\beta~=~ {-\nu_0\over \sqrt{\nu_1^2 + \nu_2^2}}\,.\eeq
This leads us to consider the scalar measure
\bel{mu1}
\mu~\doteq~ \sqrt{\nu_1^2+\nu_2^2} \cdot E\left({-\nu_0\over \sqrt{\nu_1^2 + \nu_2^2}}\right) \cdot \H^2\rest\partial^*\Omega\,,\eeq
and its projection $\Tilde\mu$ on the $t$-axis, defined by
\bel{mu2}
\Tilde \mu(S)~=~\mu\Big(\bigl\{ (t,x)\,;~t\in S, ~x\in  \R^2\bigr\}\Big)\eeq
for every Borel set $S\subset \,]0,T[\,$.   We can now define a cost functional $\Psi (\Omega)$,
for every $\Omega\in\A$, by setting
\bel{PU}\Psi (\Omega)~\doteq~\left\{\bega{cl} 
\ds \int_0^T \phi\bigl(\E(t)\bigr)\, dt \quad &  \hbox{if} ~~\Tilde\mu~\hbox{ is absolutely continuous 
with density $\E(t)$}\\
&\qquad\qquad \qquad \hbox{ w.r.t.~Lebesgue measure on $[0,T]$,}
\\[4mm]
+\infty 
&\hbox{if} ~~\Tilde\mu~\hbox{ is not absolutely continuous w.r.t.~Lebesgue measure.}
\enda\right.\eeq

\begin{remark}
{\rm In the smooth case, (\ref{PU}) yields
\bel{ceff}\Psi (\Omega)~=~\int_0^T \phi\left(\int_{\partial\Omega(t)} E(\beta)\, d\sigma\right)dt.\eeq
where $d\sigma$ denotes the 1-dimensional measure along the boundary  
$\partial \Omega(t)$.
%
%
Notice that, in the second alternative of (\ref{PU}), 
the infinite cost  is motivated by the assumption 
of superlinear growth in (\ref{phia}).
}
\end{remark}

\begin{figure}[ht]
\centerline{\hbox{\includegraphics[width=15cm]{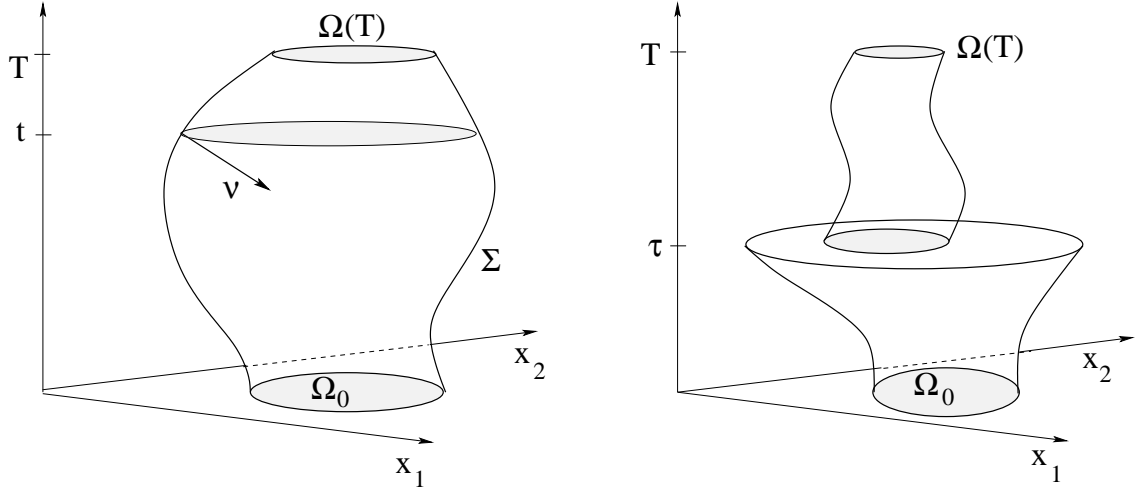}}}
\caption{\small Left: the measure $\mu$ is absolutely continuous w.r.t.~2-dimensional Hausdorff measure on the surface $\Sigma\subset\R^3$ where $u$ has a jump.
Right: an example where the projection of $\mu$ on the $t$-axis is not absolutely continuous
w.r.t.~1-dimensional Lebesgue measure.    Here $\Tilde \mu$ has a point mass at $t=\tau$.}
\label{f:sc49}
\end{figure}

\begin{remark}  {\rm If the assumptions {\bf (A1)}   hold,  the function 
\bel{LL} L(\nu)~=~ \sqrt{\nu_1^2+\nu_2^2} \cdot E\left( {- \nu_0\over \sqrt{\nu_1^2+\nu_2^2}}
\right)\eeq
is uniformly bounded on the set of all unit vectors $\nu= (\nu_0, \nu_1, \nu_2)\in \R^3$.  Hence the measure $\mu$, introduced at (\ref{mu1}), is 
absolutely continuous w.r.t.~2-dimensional Hausdorff measure on the reduced boundary
$\partial^*\Omega$.   However, its projection $\Tilde \mu$  on the time axis may 
not be absolutely continuous w.r.t.~1-dimensional Lebesgue measure, as shown in 
Fig.~\ref{f:sc49}, right.}
\end{remark}

We can now extend the function $L$ in (\ref{LL}) to a positively homogeneous function
defined for all vectors $\bfv= (v_0, v_1, v_2)\in \R^3 $, by setting
\bel{Lext}L(\bfv)~\doteq~|\bfv| \, L\left( {\bfv\over |\bfv|}\right)~=~  \sqrt{v_1^2+v_2^2} \cdot E\left( {- v_0\over \sqrt{v_1^2+v_2^2}}
\right).\eeq
To analyze the lower semicontinuity of the corresponding integral functional, we 
check whether $L$ is convex.   Restricted to the set where 
$${-v_0\over \sqrt{v_1^2+v_2^2} }~ >~ \beta_0\,,$$
the partial derivatives are
$$L_{v_0}~=~  \sqrt{v_1^2+v_2^2} \cdot E'\left( {- v_0\over \sqrt{v_1^2+v_2^2}}\right)
\cdot  {- 1\over \sqrt{v_1^2+v_2^2}}~=~ - E'\left( {- v_0\over \sqrt{v_1^2+v_2^2}}\right),$$
$$L_{v_i}~=~{v_i\over \sqrt{v_1^2+v_2^2}} \cdot E\left( {- v_0\over \sqrt{v_1^2+v_2^2}}\right)
+ {v_0 v_i\over v_1^2 + v_2^2}  E'\left( {- v_0\over \sqrt{v_1^2+v_2^2}}\right),
\qquad\qquad i=1,2,$$
$$L_{v_0v_0}~=~  E''\left( {- v_0\over \sqrt{v_1^2+v_2^2}}\right)\cdot  {1\over \sqrt{v_1^2+v_2^2}}~\geq~0\,,$$
$$L_{v_0v_i}~=~  E''\left( {- v_0\over \sqrt{v_1^2+v_2^2}}\right)\cdot (- v_0 v_i)\cdot (v_1^2+v_2^2)^{-3/2}\,,$$
$$\bega{rl}  L_{v_iv_i}&=\ds{1\over\sqrt{v_1^2+v_2^2}}\Big(1-{v_i^2\over v_1^2+v_2^2 }\Big)\,E\left( {- v_0\over \sqrt{v_1^2+v_2^2}}\right) + {v_0(v_1^2 + v_2^2) - v_0 v_i^2\over (v_1^2+v_2^2)^2}  E'\left( {- v_0\over \sqrt{v_1^2+v_2^2}}\right)
\\[3mm]&\qquad \ds+{v_0^2 v_i^2 \over (v_1^2+v_2^2)^{5/2}}E''\left( {- v_0\over \sqrt{v_1^2+v_2^2}}\right),\enda$$
$$ \bega{lr} \ds L_{v_iv_j}= -{v_iv_j\over (v_1^2+v_2^2)^{3/2}}E\left( {- v_0\over \sqrt{v_1^2+v_2^2}}\right)-{v_0v_iv_j\over (v_1^2+v_2^2)^2} E'\left( {- v_0\over \sqrt{v_1^2+v_2^2}}\right)\\[3mm]
\qquad \qquad\ds +~{v_0^2 v_iv_j\over (v_1^2+v_2^2)^{5/2}}E''\left( {- v_0\over \sqrt{v_1^2+v_2^2}}\right), \qquad \qquad i=1,2 \quad i\neq j\enda $$
\v
Since $L$ is positively homogeneous, and 
invariant under rotations in the $v_1, v_2$ coordinates, it suffices to 
compute the Hessian matrix of partial derivatives
in the special case where $\bfv=(v_0, v_1,v_2)=(v_0, 1, 0)$.    At this particular point, 
the above computations yield
\bel{Hess}D^2L(\bfv)~=~\begin{pmatrix}
E''(-v_0) & -v_0E''(-v_0)& 0\\
-v_0E''(-v_0) &  v_0^2E''(-v_0)& 0\\
0 & 0 &v_0 E'(-v_0) + E(-v_0)
\end{pmatrix}.\eeq
The eigenvalues of this matrix are
found to be 
$$0\,,\qquad (1+v_0^2)\,  E''(-v_0)\,,\qquad v_0 E'(-v_0) +E(-v_0).$$  All of these are non-negative, provided that the effort function satisfies (\ref{Eass}). 
We notice that the vector $\bfv=(v_0, 1,0)$ is itself an eigenvector of the Hessian matrix 
$D^2L(\bfv)$ at (\ref{Hess}), with zero as corresponding eigenvalue.  This is consistent with the fact that, by (\ref{Lext}), 
$L$ is linear homogeneous along each ray through the origin.
%
%

The above computations show that
the Hessian matrix of $L$ is non-negative definite at every 
point outside the cone 
\bel{Gamma}\Gamma~=~\Big\{ (v_0, v_1, v_2)\,;~v_0\geq  - \beta_0 \sqrt{ v_1^2+v_2^2}\,
\Big\}.\eeq
By (\ref{Eprop}),   $L$ vanishes on the convex cone $\Gamma$.
We thus conclude that $L$ is convex on the entire space $\R^3$.
\v
In view of (\ref{Eprop}),  the cost function $L$ in (\ref{LL}) does not pose any restriction
on how fast the set $ \Omega(t)$ expands, but it penalizes the rate at 
which $\Omega(t)$  shrinks.  This leads to the following one-sided H\"older estimate:

\begin{lemma}\label{l:21}  Let the assumptions {\bf (A1)-(A2)} hold, and let $\Omega\in\A$ be 
such that $\Psi(\Omega)<+\infty$.
Then there exists a constant $C$ such that 
\bel{1sh}
m_2\bigl(\Omega(\tau+)\setminus\Omega(\tau'-)\bigr)~\leq~C (\tau'-\tau)^{p-1\over p}
\qquad\forall 0< \tau<\tau' <T.\eeq
\end{lemma}

{\bf Proof.} {\bf 1.}
Given $\tau< \tau'$, consider the set 
\bel{Sm}
\Sigma^-~\doteq~\Big\{ (t,x)\in \partial^*\Sigma\,;~\tau<t<\tau',~\nu_0 <0\Big\}.\eeq
Observing that $E(\beta)$ is uniformly positive for $\beta\geq 0$, define the constant
$$c_0~\doteq~\min_{|\nu|=1,\, \nu_0\geq 0}~L(\nu)~=~\min_{|\nu|=1,\, \nu_0\geq 0}
\bigg\{ 
\sqrt{\nu_1^2+\nu_2^2} \cdot E\left( {- \nu_0\over \sqrt{\nu_1^2+\nu_2^2}}\right)\bigg\}.$$
One now has the estimate
\bel{mest}m_2\bigl(\Omega(\tau+)\setminus\Omega(\tau'-)\bigr)~\leq~\H^2(\Sigma^-)~
~\leq~
{1\over c_0} \int_{\Sigma^-} L(\nu)\, d\H^2\,.
\eeq
\v
{\bf 2.}
In view of the assumption (\ref{phia}), 
using H\"older's inequality with ${1\over p} + {1\over q}=1$, we obtain
$$\bega{l}\ds 
\int_{\tau}^{\tau'} 1\cdot \E(t)\, dt~\leq~ (\tau'-\tau)^{1/q}\cdot \left(\int_{\tau}^{\tau'} \bigl[\E(t)\bigr]^p\, dt
\right)^{1/p}\\[4mm]
\ds\qquad\leq~(\tau'-\tau)^{1/q}\cdot \left(\int_{\tau}^{\tau'} \left[{1\over C_1} \phi\bigl(\E(t)\bigr) + C_2\right]\, dt
\right)^{1/p}~\leq~C_3\, (\tau'-\tau)^{1/q},
\enda$$
for a suitable constant $C_3$.   Together with (\ref{mest}), this yields
$$m_2\bigl(\Omega(\tau+)\setminus\Omega(\tau'-)\bigr)~\leq~{1\over c_0} \int_{\tau}^{\tau'} \E(t)\, dt~\leq~{C_3\over c_0}\, (\tau'-\tau)^{1/q},
$$
completing the proof.
 \endproof

\begin{remark} {\rm  If $\phi(s)= +\infty $ for $s>M$, then (\ref{1sh}) can be replaced by the 
one-sided Lipschitz estimate:
\bel{1slip}
m_2\bigl(\Omega(\tau+)\setminus\Omega(\tau'-)\bigr)~\leq~C (\tau'-\tau)
\qquad\forall 0< \tau<\tau' <T.\eeq
}
\end{remark}

\section{Existence of eradication strategies}
\label{s:3}
\setcounter{equation}{0}
The next result provides a simple condition for the solvability of the Null Controllability Problem
{\bf (NCP)}.
Thinking of $\Omega(t)$ as the region contaminated by an invasive biological species, 
this yields a strategy that eradicates the contamination in finite time.

\begin{theorem}\label{t:31} Let the assumptions {\bf (A1)} hold.
Let $\Omega_0$ be a compact set whose convex closure $\Hat \Omega_0=co \,\Omega_0$ 
has perimeter which satisfies
\bel{pco}E(0)\cdot m_1(\partial \Hat\Omega_0)~<~M.\eeq
Then the null controllability problem {\bf (NCP)} has a solution.
\end{theorem}

{\bf Proof.} {\bf 1.} 
By (\ref{pco}) and the continuity of $E$, there exists a speed $\beta_1>0$ small enough so that
\bel{pc1}E(\beta_1)\cdot m_1\bigl(\partial  \Hat \Omega_0\bigr)~\leq~M.\eeq
We now define the convex subsets
\bel{omc1} \Hat \Omega(t)~=~\Big\{ x\in \R^2\,;~~B(x, \,\beta_1 t)\subseteq 
 \Hat \Omega_0\Big\}.\eeq
Notice that these sets are obtained starting from $\Hat \Omega_0$, and letting 
every boundary point $x\in \partial  \Hat\Omega(t)$ move with inward normal speed $\beta_1$.
Since the  boundaries of these sets have  decreasing length,
the total effort required by this strategy at time $t\geq 0$ is 
$$\Hat \E(t)~=~E(\beta_1)\cdot m_1\bigl(\partial \Hat\Omega(t)\bigr)
~\leq~E(\beta_1)\cdot m_1\bigl(\partial  \Hat\Omega_0\bigr)~\leq~M.$$
The set $ \Hat\Omega(t)$ shrinks to
the empty set within a finite time $ T$,  which can be estimated in terms of the diameter
of $\Omega_0$.
Namely, 
$$T~<~{\hbox{diam} ( \Hat\Omega_0)\over \beta_1}~=~{\hbox{diam} (\Omega_0)\over \beta_1}\,.$$
\v
{\bf 2.}
Next, consider the smaller sets
$$\Omega(t)~\doteq~ \Hat\Omega(t)\cap B\bigl(\Omega_0\,,~|\beta_0| t\bigr).$$
By construction, at each time $t\in [0,T]$ the boundary $\partial\Omega(t)$ either
touches the boundary $\partial \Hat\Omega(t)$, ore else it expands with normal speed
$|\beta_0|$.   This implies
$$\E(t)~=~\int_{\partial \Omega(t)} E\bigl(\beta(t,x)\bigr)\, d\sigma~=~\int_{\partial \Omega(t)
\cap \partial \Hat\Omega(t)} E(\beta_1)\, d\sigma~\leq~ \Hat\E(t)~\leq~M.$$
Therefore, the multifunction $t\mapsto\Omega(t)$ satisfies  (\ref{eb}).
\v
{\bf 3.} It remains to prove that the initial condition is satisfied, in the sense that 
\bel{icf}
\lim_{t\to 0+}\, \bigl\| {\bf 1}_{\Omega(t)} - {\bf 1}_{\Omega_0}\bigr\|_{\L^1(\R^2)}~=~0.\eeq
Toward this goal, we observe that, since $\Omega_0$ is compact,
\bel{ic5}
\lim_{t\to 0+} m_2\Big( B\bigl(\Omega_0,\, 2|\beta_0| t\bigr)\Big) ~=~m_2(\Omega_0).\eeq
Hence
\bel{ge1}\lim_{t\to 0+} m_2\bigl( \Omega(t)\setminus\Omega_0\bigl) ~=~ 0.\eeq
On the other hand, observing that
$$\Omega_0\setminus \Omega(t)~\subseteq~\bigl\{
x\in \Hat\Omega_0\,;~~d(x, \partial \Hat\Omega_0)~\leq~|\beta_1|t\Bigr\},$$
we obtain
\bel{ge2}
\lim_{t\to 0+} m_2\bigl( \Omega_0\setminus \Omega(t)\bigr) ~\leq~
\lim_{t\to 0+} |\beta_1|\,t\cdot m_1(\partial \Hat\Omega)~=~ 0.\eeq
Together, (\ref{ge1}) and (\ref{ge2}) yield (\ref{icf}).
\endproof

\section{Existence of optimal strategies}
\label{s:4}
\setcounter{equation}{0}

To achieve the existence of an optimal strategy $t\mapsto\Omega(t)$, we 
need to somewhat relax the formulation of the problem
{\bf (OP)}. 
We recall that a subset $\Omega\subset \,]0,T[\,\times \R^2$ determines
the multifunction $t\mapsto \Omega(t)$ as in (\ref{slice}).  Moreover,
assuming that $\Omega\in\A$ is bounded and has finite perimeter,   the initial and terminal 
values $\Omega(0)$
and $\Omega(T)$ are uniquely determined by (\ref{om0T}).
Recalling the functional $\Psi(\Omega)$ introduced at (\ref{PU}), and denoting by
$m_2$, $m_3$ respectively the 2- and 3-dimensional Lebesgue measure,
we thus consider the problem of Optimal Set Motion:
\begi
\item[{\bf (OSM)}] {\it Given a bounded initial set $\Omega_0\subset\R^2$,
find a set  $\Omega\subset \,]0,T[\,\times \R^2$
which minimizes the functional
\bel{FU}
\J (\Omega)~\doteq~\Psi(\Omega) + c_1 m_3(\Omega) + c_2 \, m_2\bigl(\Omega(T)\bigr)\,,\eeq
among all sets $\Omega\in \A$ such that 
$\Omega(0)=\Omega_0$.}
 \endi
Aim of this section is to prove the existence of solutions to the above optimization 
problem.

\begin{theorem}\label{t:41}
Let  $E, \phi$ satisfy the assumptions {\bf (A1)-(A2)}.
Then, for any compact set $\Omega_0\subset\R^2$ with finite perimeter and any $T, c_1,c_2>0$,
the problem {\bf (OSM)} has an optimal solution.
\end{theorem}

{\bf Proof.} {\bf 1.}  We start with the trivial observation that $\J(\Omega)\geq 0$ for every 
$\Omega\in \A$.
Moreover, choosing
$\Omega = \,]0,T[\,\times\Omega_0$, so that 
$\Omega(t)\equiv{\Omega_0}$ for all $t\in [0,1]$,  we obtain an admissible 
set $\Omega \in \A$
with $\J(\Omega)< +\infty$.    We can thus consider 
a minimizing sequence of sets  $\Omega_n\in \A$ 
such that, as $n\to\infty$,
$$\J(\Omega_n)~\to~\J_{min}~\doteq~\inf_{\Omega\in \A} \J (\Omega)\,.$$
Without loss of generality we can assume that the sets 
$$\Omega_n(t)~\doteq~\bigl\{x\in \R^2\,; ~(t,x)\in\Omega_n\bigr\}$$
are contained in the neighborhood of radius $|\beta_0| t$ around $\Omega_0$:
\bel{omn}
\Omega_n(t)~\subseteq~B\bigl( \Omega_0\,;~|\beta_0| t\bigr)\eeq
for every $n\geq 1$ and $t\geq 0$.   Otherwise, we can simply replace 
each set $\Omega_n(t)$ with the intersection
$\Omega_n(t)\cap B\bigl( \Omega_0\,;~|\beta_0| t\bigr)$, without increasing 
the total cost.
\v
{\bf 2.} In the next two steps we prove a uniform bound on perimeters of the 
sets $\Omega_n\subset \R^3$.  

\begin{figure}[ht]
\centerline{\hbox{\includegraphics[width=13cm]{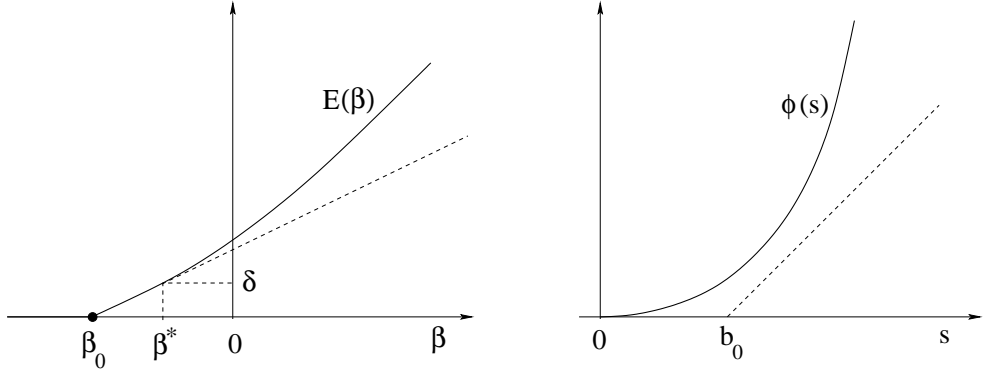}}}
\caption{\small Left: a lower bound on the function $E(\beta)$, considered at (\ref{Eb}).
  Right: since the cost function $\phi$ has superlinear growth, it admits a lower bound
of the form (\ref{phs}).}
\label{f:sc46}
\end{figure}

Choose a speed $\beta^*<0$ and constants $\delta,\lambda>0$ such that (see Fig.~\ref{f:sc46}, left)
\bel{Eb}
E(\beta^*)~=~\delta ~>~0,\qquad 
E(\beta)~\geq~\delta +\lambda\,(\beta-\beta^*)\qquad\forall \beta\geq \beta^*\,.
\eeq
We split the reduced boundary  in the form
\bel{spn}
\partial^*\Omega_n~=~\Sigma_n^-\cup\Sigma_n^+\,,\eeq
so that the following holds.
Calling $\nu=(\nu_0, \nu_1,\nu_2)$ the normal vector at the point
$(t,x)\in \partial^*\Omega$, and defining 
 the inner normal velocity $\beta_n=\beta_n(t,x)$ as in 
(\ref{bdef}), one has
\bel{Snpm}
\left\{ \bega{rl} \beta_n(t,x)~\leq ~\beta^*\qquad &\hbox{if}\quad (t,x)\in \Sigma_n^-\,
,\\[2mm] \beta_n(t,x)~>~\beta^*\qquad &\hbox{if}\quad (t,x)\in \Sigma_n^+\,.
\enda\right.\eeq

By {\bf (A2)} we can find a constant $b_0>0$ such that
\bel{phs}
\phi(s)~\geq~s-b_0\qquad\qquad\forall s\in\R_+\,.\eeq
In view of (\ref{Eprop}), this implies
\bel{Etb}\bega{rl}\ds
\int_0^T\phi\bigl(\E_n(t)\bigr)\, dt&\ds \geq~\int_{\Sigma_n} \sqrt{\nu_1^2+\nu_2^2} \cdot E\left({-\nu_0\over \sqrt{\nu_1^2 + \nu_2^2}}\right) \, d\H^2 - b_0 T.
\enda  \eeq
On the domain $\Sigma_n^+$, where 
\bel{bbig}\beta~\doteq~{-\nu_0\over \sqrt{\nu_1^2+\nu_2^2}}~\geq~\beta^*,\eeq
by (\ref{Eb})  we have
\bel{st1}
 \sqrt{\nu_1^2+\nu_2^2} \cdot E\left({-\nu_0\over \sqrt{\nu_1^2 + \nu_2^2}}\right)
 ~\geq~\sqrt{\nu_1^2+\nu_2^2}\left( \delta + \lambda\Big( {-\nu_0\over \sqrt{\nu_1^2 + \nu_2^2}}- \beta^*\Big)\right)~\geq~c_3
\eeq
for some constant $c_3>0$. Together with (\ref{Etb}), this yields
\bel{Sp1} \int_{\Sigma_n^+} d\H^2~\leq~{1\over c_3}\left[ 
\int_0^T\phi\bigl(\E_n(t)\bigr)\, dt + b_0 T\right] .\eeq

\v
{\bf 3.}
On the domain $\Sigma_n^-$ where (\ref{bbig}) fails,
we have a lower bound
\bel{n00}
\nu_0~\geq~c_4~>~0,\eeq
for some positive constant $c_4$.  We can now write
\bel{Sm1}\bega{l}\ds m_2\Big( B\bigl(\Omega_0,\,|\beta_0| T\bigr)\Big)~\geq~
m_2\bigl(\Omega_n(T)\bigr) - m_2\bigl(\Omega_n(0)\bigr) \\[4mm]
 \qquad \ds
=~\int_{\Sigma_n^+\cup\Sigma_n^-} \nu_0\,d\H^2
~ \geq~c_4 \int_{\Sigma_n^-}
d\H^2 - \int_{\Sigma_n^+} d\H^2.
\enda
\eeq
Combining (\ref{Sm1}) with (\ref{Sp1}) one obtains
\bel{Sm2} \int_{\Sigma_n^-}
d\Sigma~\leq~{1\over c_4}\left\{ m_2\Big( B\bigl(\Omega_0,\,|\beta_0| T\bigr)\Big)
+ {1\over c_3}\left[ 
\int_0^T\phi\bigl(\E_n(t)\bigr)\, dt + b_0 T\right] \right\}.
\eeq
We notice that, for a minimizing sequence, the integrals $\int\E_n(t)dt$ 
must be bounded, because they are part of the cost
functional (\ref{PU}).  

Together, the two inequalities (\ref{Sp1}) and (\ref{Sm2}) thus yield a
uniform bound on the 2-dimensional measure $\H^2(\Sigma_n)$, i.e.~on the total variation
of the function $u_n$, for every $n\geq 1$.

\v
{\bf 4.} 
Thanks to the uniform BV bound, by possibly taking a subsequence,
a compactness argument (see Theorem~12.26 in \cite{M}) yields the the existence of a 
bounded set with finite perimeter $\Omega\in \A$ such that the following holds.
As $n\to\infty$, one has the convergence 
\bel{cvg}\Big\| {\bf 1}_{\Omega_n} - {\bf 1}_\Omega\Big\|_{\L^1\bigl(]0,T[\,\times\R^2\bigr)}~\to~0,\eeq
together with the weak convergence of measures
\bel{wco}
\mu_{\Omega_n}~\wto~\mu_\Omega\,.\eeq

Next, the assumption of superlinear growth 
(\ref{phia}) implies
\bel{L1b}\int_0^T\phi\bigl(\E_n(t)\bigr)dt~ \geq~ C_1\int_0^T\bigl[\E_n(t)
\bigr]^p\,dt - C_2 T.\eeq
Therefore, the functions $\E_n$ are uniformly bounded in $\L^p$. Since $p>1$,
we can extract a weakly convergent subsequence 
$\E_n\wto \E_\infty\in\L^p([0,T],\R)$. 
We observe that  $\E_\infty$ yields the density of the measure 
$\Tilde\mu_\infty$, defined as the weak limit of the projected measures $\Tilde\mu_n$. 
\v
{\bf 5.}
To prove a lower semicontinuity result, we first replace $\phi$ with a globally Lipschitz, convex function
\bel{pms}
\phi_m(s)~=~\sup_{0\leq a\leq m,\, b\in\R}   \Big\{ as+b\,;~~at+b\leq \phi(t)\quad\forall t\in\R\Big\}.\eeq
Notice that the function 
$\phi_m$ is Lipschitz continuous with constant $m$. Its graph is the upper envelope of all straight lines with slope $\leq m$ that lie below the graph of $\phi$.

Let $\ve>0$ be given. Choosing $m$ large enough, we achieve
\bel{ap1}
\int_0^T \phi\bigl(\E_\infty(t)\bigr)\, dt~\leq~\ve +
\int_0^T \phi_m\bigl(\E_\infty(t)\bigr)\, dt .\eeq
\v
{\bf 6.} Next, we partition  the interval $[0,T]$ 
into finitely many subintervals $I_j = [t_{j-1}, t_j]$, $j=1,\ldots, N$, so that the difference
between $\E_\infty$ and its average value over each subinterval is bounded by
\bel{ap2}
\sum_{1\leq j\leq N} \int_{t_{j-1}}^{t_j} \left| \E_\infty(t) - \avint_{t_{j-1}}^{t_j}  \E_\infty(\tau)d\tau\right| dt
~<~{\ve\over m}\,.\eeq
The Lipschitz continuity of the function $\phi_m$ now yields
\bel{ap3}
\sum_{1\leq j\leq N} \int_{t_{j-1}}^{t_j} \left| \phi_m\bigl(\E_\infty(t)\bigr) - \phi_m\left(\avint_{t_{j-1}}^{t_j}  \E_\infty(\tau)d\tau\right)\right| dt
~<~\ve\,,\eeq
\bel{ap4}
\int_0^T \phi_m  \bigl(\E_\infty(t)\bigr)\,dt~\leq
~\ve+ \sum_{1\leq j\leq N}(t_j-t_{j-1}) 
 \phi_m\left(
 \avint_{t_{j-1}}^{t_j}  \E_\infty(\tau)d\tau\right).
\eeq
\v
{\bf 7.} We now study the relation between $\E_\infty(t)$ and the instantaneous effort $\E(t)$ associated to the 
limit set $\Omega$.   

Since the function 
$L$ in (\ref{LL}) is convex, for any $0<\tau<\tau'<T$
we can use a lower semicontinuity result for anisotropic functionals
 (see Theorem 20.1 in \cite{M})  and conclude 
\bel{lsc}\bega{rl} \ds
\int_\tau^{\tau'}  \E(t)\, dt &\ds
=~\int_{\partial^*\Omega\cap\{ \tau<t<\tau'\}}   L(\nu)\, d\H^2~\leq~\liminf_{n\to\infty} 
\int_{\partial^*\Omega\cap\{ \tau<t<\tau'\}}   L(\nu_n)\, d\H^2\\[4mm]
&\ds =~\liminf_{n\to\infty}  \int_\tau^{\tau'}  \E_n(t)\, dt~=~ \int_\tau^{\tau'} 
 \E_\infty(t)\, dt\,.\enda \eeq
Since $\tau, \tau'$ were arbitrary, this implies $\E(t)\leq \E_\infty(t)$ for a.e.~$t\in [0,T]$.
\v
Next, 
 by Jensen's inequality and the convexity of $\phi_m$  it follows 
\bel{ap6}(t_j-t_{j-1}) 
 \phi_m\left(
 \avint_{t_{j-1}}^{t_j}  \E_n(\tau)d\tau\right)~\leq~\int_{t_{j-1}}^{t_j} 
\phi_m\bigl( \E_n(t)\bigr)d t\,.\eeq
Summing over $j=1,\ldots,N$, and using (\ref{ap4}) and (\ref{ap1}), (\ref{ap3}), we conclude
\bel{limp}\bega{l} \ds
\int_0^T \phi\bigl(\E(t)\bigr)\,dt~ \leq~\int_0^T \phi\bigl(\E_\infty(t)\bigr)\,dt~\leq~
\ve + \int_0^T \phi_m\bigl(\E_\infty(t)\bigr)\,dt
\\[4mm]
\ds \qquad \leq 2\ve+ \sum_{1\leq j\leq N}(t_j-t_{j-1}) 
 \phi_m\left(
 \avint_{t_{j-1}}^{t_j}  \E_\infty(\tau)d\tau\right)
 \\[4mm]
 \ds \qquad = 2\ve+ \sum_{1\leq j\leq N}(t_j-t_{j-1}) 
 \phi_m\left(\lim_{n\to\infty}
 \avint_{t_{j-1}}^{t_j}  \E_n(\tau)d\tau\right)\\[4mm]
\qquad \ds =~2\ve +\lim_{n\to\infty} \sum_{1\leq j\leq N}(t_j-t_{j-1}) 
 \phi_m\left(
 \avint_{t_{j-1}}^{t_j}  \E_n(\tau)d\tau\right)\\[4mm]
 \qquad \ds \leq~2\ve +\liminf_{n\to \infty} \int_0^T \phi_m\bigl(\E_n(t)\bigr)\,dt
 ~\leq~
2\ve + \liminf_{n\to \infty} \int_0^T \phi\bigl(\E_n(t)\bigr)\,dt\,.
\enda\eeq
\v
{\bf 8.} As $n\to\infty$, the convergence (\ref{cvg}) immediately implies
\bel{cm3}m_3(\Omega_n)~\to~m_3(\Omega).\eeq
    Moreover, since the map $t\mapsto {\bf 1}_{\Omega(t)}\in \L^1(\R^2)$ has bounded variation, given $\ve>0$ we can find $\delta>0$ such that 
\bel{m2o}
\left| m_2(\Omega(T)\big)-  {1\over\delta} \int_{T-\delta}^T m_2(\Omega(t)\big)\, dt
\right|~<~\ve.\eeq
On the other hand, by the one-sided estimate in Lemma~\ref{l:21} it follows
$$m_2\bigl(\Omega_n(\tau)\setminus\Omega_n(T)\bigr)~\leq~ C\,(T-\tau)^{p-1\over p},
$$
for some constant $C$ independent of $T,\tau$, and $n$.  Therefore, we can find $\delta>0$ 
such that
$$m_2\bigl(\Omega_n(T)\bigr)~\geq~{1\over \delta} \int_{T-\delta}^T m_2
\bigl(\Omega_n(t)\bigr)\, dt -\ve$$
for every $n\geq 1$.  Since
$${1\over \delta} \int_{T-\delta}^T m_2
\bigl(\Omega_n(t)\bigr)\, dt ~\to~ {1\over\delta} \int_{T-\delta}^T m_2
\bigl(\Omega(t)\bigr)\, dt \,,$$
we conclude
\bel{lsom}
m_2\bigl(\Omega(T)\bigr) ~\leq~\liminf_{n\to\infty}  m_2\bigl(\Omega_n(T)\bigr)  + 2\ve.\eeq
\v
{\bf 9.} Combining (\ref{limp}), (\ref{cm3}), and (\ref{lsom}),
since $\ve>0$ was arbitrary we conclude 
\bel{linf}
\J(\Omega)~\doteq~
\int_0^T\E(t)\, dt + c_1\, m_3(\Omega) + c_2 \, m_2\bigl(\Omega(T)\bigr) ~\leq~\liminf_{n\to \infty} \J(\Omega_n).\eeq
\v
{\bf 10.} It remains to prove that the limit set $\Omega$ satisfies the initial 
condition (\ref{icf}).   Since $\Omega_0$ is compact, by (\ref{omn}) we immediately have
\bel{li4}\lim_{t\to 0+} m_2\bigl(\Omega(t)\setminus\Omega_0\bigr)~\leq~\lim_{t\to 0+} 
m_2\Big(B\bigl(\Omega_0;\, |\beta_0|t\bigr) \setminus\Omega_0\Big) ~=~0.\eeq
On the other hand,   by (\ref{1sh}) for every $t>0$ one has
$$m_2\bigl( \Omega_0\setminus \Omega_n(t)\bigr)  ~\leq~C\, t^{p-1\over p},
$$
for a suitable constant $C$  independent of $t$ and $ n$.
Taking the limit as $n\to\infty$ one obtains
\bel{li9} 
m_2\bigl(\Omega_0\setminus \Omega(t)\bigr)~ \leq~C\, t^{p-1\over p}.\eeq
Together, (\ref{li4}) and (\ref{li9}) yield the convergence ${\bf 1}_{\Omega(t)}\to 
{\bf 1}_{\Omega_0}$ in $\L^1(\R^2)$, completing the proof.
\endproof
%
%

By entirely similar arguments one can prove the existence of an optimal solution for the minimum time problem.

\begin{theorem}\label{t:42}
Let the functions $E$ satisfy the assumptions {\bf (A1)} and let $M>0$ be given.
Let $\Omega_0\subset\R^2$ be a compact set with finite perimeter such that 
the null controllability problem {\bf (NCP)} has a solution.   Then the minimum time problem
{\bf (MTP)} has an optimal solution.
\end{theorem}

{\bf Proof.} Consider a minimizing sequence $(\Omega_n)_{n\geq 1}$.   Calling 
$\Omega_n(t)=\bigl\{x\,;~(t,x)\in\Omega_n\bigl\}$, we thus have
\bel{mom}
m_2\bigl(\Omega_n(T_n)\bigr)~=~0,\eeq
with $T_n\downarrow T$.     The same arguments used in the proof of Theorem~\ref{t:41}
yield a uniform bound on the perimeter of $\Omega_n$.   By possibly taking a subsequence
we obtain the strong 
convergence ${\bf 1}_{\Omega_n}\to {\bf 1}_\Omega$.   By assumption, $\E_n(t)\leq M$ for every $t,n$.   Calling 
$\E(t)$ the effort corresponding to $\Omega$, and 
$\E_\infty$ the weak limit of the functions $\E_n$, the previous arguments
yield $\E(t)\leq \E_\infty(t)\leq M$ for all $t\geq 0$.   

In the present setting,  for any $n\geq 1$, we have the 
one-sided Lipschitz estimate (\ref{1slip}).
In particular, taking  $\tau=T_n$ and $\tau'=T$,  we obtain
$$m_2\big(\Omega_n(T) \bigr)~\leq~C_4(T_n-T).$$
Taking the limit as $n\to\infty$, this implies $m_2\bigl(\Omega(T)\bigr)=0$.  
Hence the set-valued function $t\mapsto \Omega(t) $ provides an optimal solution
to the minimum time problem.
\endproof

\section{Necessary conditions for optimality}
\label{s:5}
\setcounter{equation}{0}

\begin{figure}[ht]
\centerline{\hbox{\includegraphics[width=10cm]{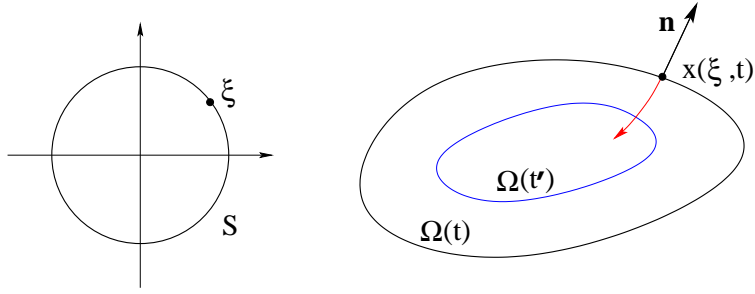}}}
\caption{\small   At each time $t\in [0,T]$, the boundary of the set $\Omega(t)$
is parameterized by $\xi\mapsto x(t,\xi)$, where $\xi\in S$ ranges over the unit circle.}
\label{f:sc6}
\end{figure}

Let $t\mapsto \Omega(t)\subset\R^2$ be an optimal solution for the problem 
{\bf (OP)} of control of a moving set. 
Aim of this section is to derive a set of necessary conditions for optimality, in the form of
a Pontryagin maximum principle \cite{BPi, Ce, FR}.   For this purpose, somewhat stronger regularity assumptions will be needed.
As shown in Fig.~\ref{f:sc6},
we consider the unit circumference $S=\{\xi\in \R^2\,;~|\xi|=1\}$ and, 
for each $t\in [0,T]$, we  assume that
$$\xi~\mapsto ~x(t,\xi)~\in ~\partial \Omega(t)$$ is a $\C^2$ parameterization of the 
boundary of $\Omega(t)$ (oriented counterclockwise), satisfying 
\begi
\item[{\bf (A3)}] {\it  There exists a constant $C>0$ such that
\bel{ra1} {1\over C} ~\leq~\bigl|x_\xi(t,\xi)\bigr|~\leq~C\qquad\forall (t,\xi)\in [0,T]\times S,\eeq
Moreover, for every $\xi\in S$ the trajectory
$t\mapsto x(t,\xi)$ is orthogonal to the boundary $\partial\Omega(t)$ at every time $t$.   Namely,
\bel{xt1}x_t(t,\xi)~=~\beta(t,\xi)\,\bfn(t,\xi),\eeq
where $\bfn= (n_1, n_2)$ is the unit inner normal to $\partial \Omega(t)$ at the point $x(t,\xi)$, 
and $\beta$ is a continuous, scalar function.}
\endi

Throughout the following, we write $\bfn^\perp=(- n_2, n_1)$ for the perpendicular 
vector and denote by 
\bel{curv}
\omega(t,\xi)~\doteq~{1\over \bigl|x_\xi(t,\xi)\bigr|}\, \la \bfn^\perp(t,\xi),\bfn_\xi(t,\xi)\ra\eeq
the curvature of the boundary $\partial \Omega(t)$ at the point $x(t,\xi)$.

To derive a set of optimality conditions, we introduce the
adjoint function $Y:[0,T]\times S\mapsto \R$, defined as the solution of
the linearized equation
\bel{Ydt}
Y_t(t,\xi) ~=~ \left(\beta(t,\xi) -{E(\beta(t,\xi))\over E'(\beta(t,\xi))}\right)\omega(t,\xi) \,Y(t,\xi)  - c_1 \,,
\eeq
with terminal condition
\bel{YT} Y(T,\xi)~=~c_2.\eeq
Notice that (\ref{Ydt}) yields a family of linear ODEs, that can be independently solved for
each $\xi\in S$.
In addition, we consider the function
\bel{lamb}\lambda(t)~\doteq~\phi'\bigl(\E(t)\bigr)~=~\phi'\left( \int_S  E(\beta(t,\xi)) \, |x_\xi(t,\xi)|\, d\xi\right).\eeq

We are now ready to state our main result, providing necessary conditions for optimality.

\begin{figure}[ht]
\centerline{\hbox{\includegraphics[width=9cm]{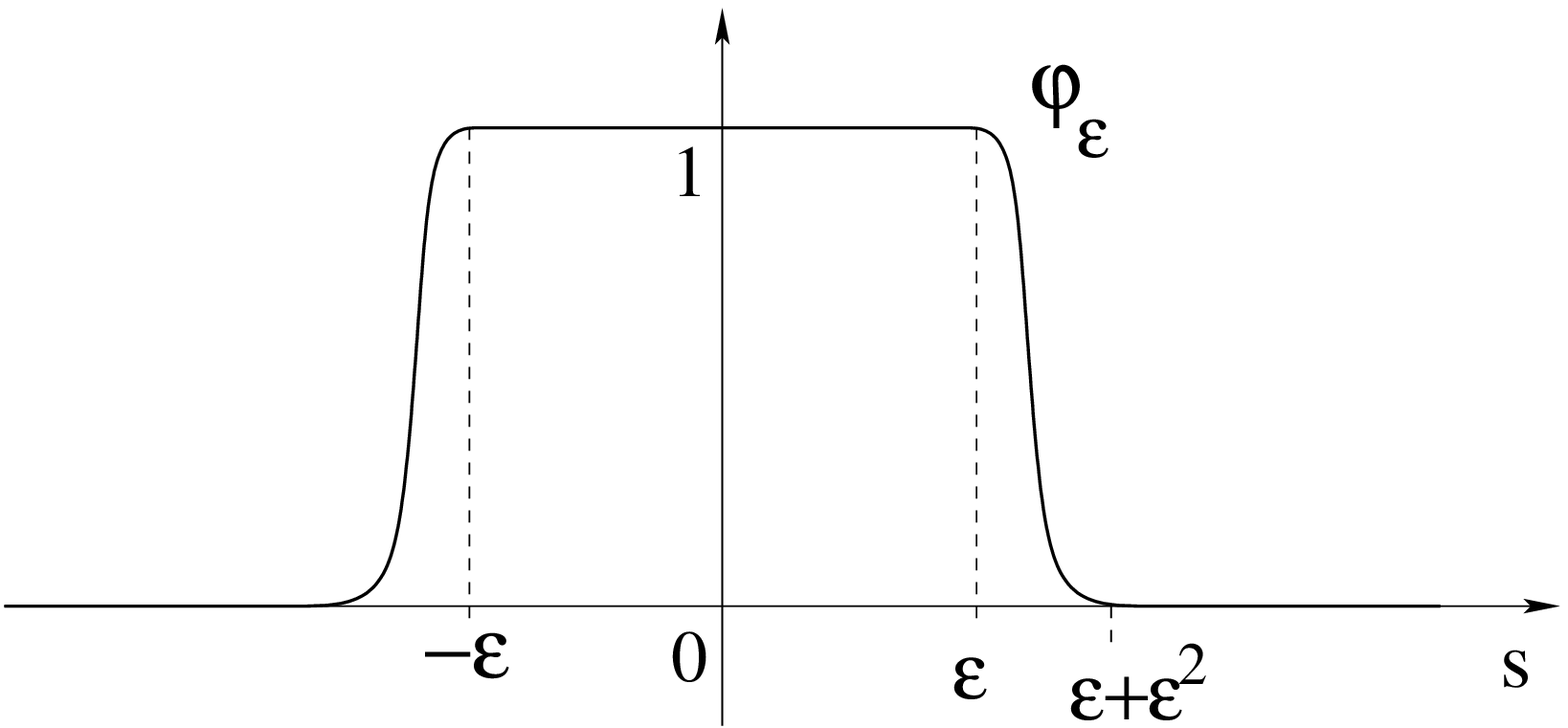}}}
\caption{\small The functions $\vp_\ve$ introduced at (\ref{vpep}). }
\label{f:sc14}
\end{figure}

\begin{figure}[ht]
\centerline{\hbox{\includegraphics[width=9cm]{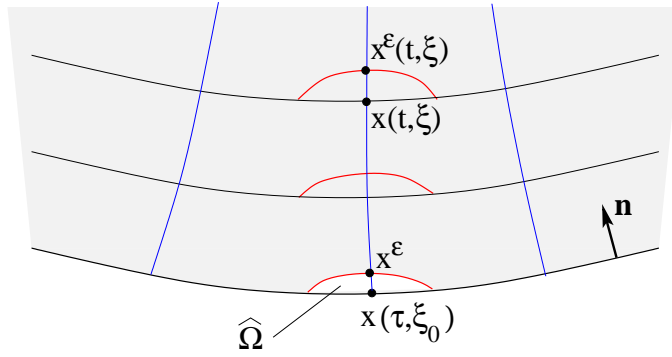}}}
\caption{\small  A perturbation of the optimal  strategy.  If at time $\tau$ the additional 
region $\Hat\Omega$ could be freed from the contamination, the total cost would be reduced in the amount $m_2(\Hat\Omega)\cdot Y(\tau, \xi_0)$. }
\label{f:sc21}
\end{figure}

\begin{theorem} \label{t:51}  Let $E$ satisfy {\bf (A1)} and let $\phi:\R_+\mapsto\R_+$ 
be a $\C^1$ function which satisfies {\bf (A2)}. 
Assume that $t\mapsto \Omega(t)$ provides an optimal solution to {\bf (OP)}.   Let $\xi\mapsto x(t,\xi)$ be a $\C^2$ parameterization of the boundary
of the set $\Omega(t)$, satisfying  {\bf (A3)}. Let $\lambda(t)$ be as in (\ref{lamb})
and consider the adjoint function $Y=Y(t,\xi)$  constructed at (\ref{Ydt})-(\ref{YT}).

Then, for every 
$t\in [0,T]$ and $\xi\in S$,  the normal velocity $\beta=\beta(t,\xi)$ satisfies
\bel{max1}
\lambda(t) E\big(\beta(t,\xi)\bigr) -Y(t,\xi) \beta(t,\xi)  ~= ~\min_{\beta\geq \beta_0}~
\Big\{\lambda(t) E(\beta)-Y(t,\xi) \beta  \Big\}.\eeq
\end{theorem}
\v
{\bf Proof.}  {\bf 1.} Assume that the conclusion fails.  
We can thus  find  a point $(\tau, \xi_0)\in ]0,T[\,\times S$ 
and some value $\ov \beta\in \R$ such that 
\bel{max2}
\lambda(\tau) E\big(\beta(\tau, \xi_0)\bigr) -Y(\tau, \xi_0) \beta(\tau, \xi_0)  ~> ~
\lambda(\tau) E(\ov \beta)-Y(\tau, \xi_0) \ov \beta  .\eeq
 We shall construct a family of perturbations $\Omega^\ve(t)$, $t\in [0,T]$, which achieve a lower cost. Here $\Omega^\ve(t)$ will be the set with boundary
\bel{POme}
\partial \Omega^\ve(t)~=~\bigl\{ x^\ve(t,\xi)\,;~\xi\in S\bigr\},\eeq
for suitable perturbations $x^\ve$  (see Fig.~\ref{f:sc21}).  These resemble the ``needle variations" used in the 
classical proof of the Pontryagin maximum principle.  Namely,  we perform a large change
in the inward normal velocity $\beta=\beta(t,\xi)$ (which here plays the role of a control function) on the small domain 
$[\tau-\ve^4, \tau]\times [\xi_0-\ve, \xi_0+\ve]$.   At all subsequent times $t\in [\tau,T]$, 
we choose
the perturbed normal velocity
$\beta^\ve$  so that its cost 
remains almost the same as is the original solution.
\v
{\bf 2.}
As a preliminary, consider a smooth function
$\vp:\R\mapsto [0,1]$ such that 
$$\left\{ \bega{rl} \vp(s)=1\quad &\hbox{if}~~ s\leq 0,\cr
\vp(s)=0\quad &\hbox{if}~~ s\geq 1,\cr
\vp'(s)\leq 0\quad &\hbox{for all}~~ s\in \R.\enda\right.$$
As shown in Fig.~\ref{f:sc14}, we then define the functions
\bel{vpep}\vp_\ve(s)~\doteq~\vp\left( {|s|-\ve\over\ve^2}\right).\eeq
Recalling (\ref{curv}),  for $t\in [\tau,T]$ we define $X=X(t,\xi)$ to be 
the solution to the linearized evolution equation
\bel{Xt2}
X_t(t,\xi)~=~{E\bigl(\beta(t,\xi)\bigr)\over 
E'\bigl(\beta(t,\xi)\bigr)}\cdot \omega(t,\xi)\, X(t,\xi),\eeq
with initial data at $t=\tau$ given by
\bel{Xt3}
X(\tau,\xi)~=~\bigl[ \ov \beta - \beta(\tau, \xi_0)\bigr]\cdot \vp_\ve(\xi-\xi_0).\eeq

 The perturbations $x^\ve$ are now defined as follows.
On an  initial time interval, we set
\bel{xe1}
x^\ve(t,\xi)~=~x(t,\xi)\qquad \hbox{if}\quad t\in [0, \tau-\ve^4], ~\xi\in S.\eeq
For $t\geq \tau$, recalling  (\ref{Xt2})-(\ref{Xt3}), we define
\bel{xe2}
x^\ve(t,\xi)~=~x(t,\xi)+ \ve^4 X(t,\xi)\bfn(t,\xi)\qquad \hbox{if}\quad t\in [\tau, T], ~\xi\in S.\eeq
Finally, in the remaining small interval of time before $\tau$, we define $x^\ve$ 
by setting
\bel{xe3}
x^\ve(t,\xi)~=~x(t,\xi)+ \bigl[ \ve^4 - (\tau-t)\bigr]\,X(\tau,\xi)\bfn(t,\xi)\qquad \hbox{if}\quad t\in [\tau-\ve^4, \,\tau], ~\xi\in S.\eeq
 In the remainder of the proof we will show that, if (\ref{max2}) holds, then 
 for a suitably small $\ve>0$ the perturbed sets
 $\Omega^\ve(t)$ in 
(\ref{POme}) achieve a strictly smaller cost.
In view of the definition of the perturbation $X$, it will be convenient to split the domain as
\bel{split1}S~=~S_1\cup S_2\cup S_3\,,\eeq
where
\bel{split2}\bega{rl} S_1&\doteq~\bigl\{ \xi\in S\,;~|\xi-\xi_0|\leq \ve\bigr\},\\[2mm]
S_2&\doteq~\bigl\{ \xi\in S\,;~\ve<|\xi-\xi_0|< \ve+\ve^2\bigr\},\\[2mm]
S_3&\doteq~\bigl\{ \xi\in S\,;~|\xi-\xi_0|\geq \ve+\ve^2\bigr\}.\enda
\eeq

%

\v
{\bf 3.} We begin by analyzing what happens during the interval $[\tau-\ve^4, \tau]$,
where, by (\ref{Xt3}) and (\ref{xe3}),
\bel{xe4}x^\ve(t,\xi)~=~x(t,\xi)
+ \bigl[ \ve^4 - (\tau-t)\bigr]\,\bigl[ \ov \beta - \beta(\tau, \xi_0)\bigr]\cdot \vp_\ve(\xi-\xi_0)\bfn(t,\xi).\eeq
Differentiating w.r.t.~$\xi$ and recalling (\ref{vpep}),  we obtain
\bel{xe5}\bega{rl}x^\ve_\xi(t,\xi)&=~x_\xi(t,\xi)
+ \bigl[ \ve^4 - (\tau-t)\bigr]\,\bigl[ \ov \beta - \beta(\tau, \xi_0)\bigr]\cdot \vp'_\ve(\xi-\xi_0)\bfn(t,\xi)\\[3mm]
&\qquad  +\bigl[ \ve^4 - (\tau-t)\bigr]\,\bigl[ \ov \beta - \beta(\tau, \xi_0)\bigr]\cdot \vp_\ve(\xi-\xi_0)\bfn_\xi(t,\xi) .\enda\eeq
In connection with the decomposition (\ref{split1})-(\ref{split2}), by
(\ref{Xt3}), we have
\bel{xe6}
x^\ve_\xi(t,\xi)- x_\xi(t,\xi)~=~\left\{\bega{cl} \O(1)\cdot \ve^4\qquad &\hbox{if}\quad \xi\in S_1\,,\\[1mm]
\O(1)\cdot \ve^2\qquad &\hbox{if}\quad\xi\in S_2\,,\\[1mm]
0\qquad &\hbox{if}\quad \xi\in S_3\,.
\enda\right.
\eeq
In view of (\ref{ra1}), the same bounds hold for the normal vector $\bfn^\ve(t,\xi)$, namely
\bel{xn}
\bfn^\ve(t,\xi)-\bfn(t,\xi)~=~\left\{\bega{cl} \O(1)\cdot \ve^4\qquad &\hbox{if}\quad \xi\in S_1\,,\\[1mm]
\O(1)\cdot \ve^2\qquad &\hbox{if}\quad\xi\in S_2\,,\\[1mm]
0\qquad &\hbox{if}\quad \xi\in S_3\,.
\enda\right.
\eeq
Next, 
differentiating (\ref{xe4}) w.r.t.~time,
we compute
\begin{equation}\label{xe77}
\bega{rl}
\beta^\ve(t,\xi)&=~\la x^\ve_t(t,\xi) \,,~\bfn^\ve(t,\xi)\ra\\[3mm]
&=~ \Big\langle \beta(t,\xi)\bfn(t,\xi) + \bigl[\ov \beta - \beta(\tau,\xi_0)\bigr] \cdot \vp_\ve(\xi-\xi_0)\bfn(t,\xi),\,\bfn^\ve(t,\xi)\Big\rangle\\[3mm]
&\qquad\qquad +~  \Big\langle  \bigl[ \ve^4 - (\tau-t)\bigr]\, \bigl[\ov \beta - \beta(\tau,\xi_0)\bigr] \cdot \vp_\ve(\xi-\xi_0)\bfn_t(t,\xi),\,\bfn^\ve(t,\xi)\Big\rangle\\[3mm]
&=\ds ~\beta(t,\xi) +  \Big\langle \beta(t,\xi)\bfn(t,\xi) ,\,\bfn^\ve(t,\xi)- \bfn(t,\xi)\Big\rangle
\\[3mm]
&\qquad\qquad \ds
+\bigl[\ov \beta - \beta(\tau,\xi_0)\bigr] \vp_\ve(\xi-\xi_0)]
\Big( 1+ \la \bfn(t,\xi),\, \bfn^\ve(t,\xi) - \bfn(t,\xi)\ra\Big)
\\[3mm]
&\ds\qquad\qquad + \bigl[ \ve^4 - (\tau-t)\bigr]\, \bigl[\ov \beta - \beta(\tau,\xi_0)\bigr] \vp_\ve(\xi-\xi_0)\cdot \la  \bfn_t(t,\xi),\,\bfn^\ve(t,\xi)- \bfn(t,\xi)\ra.
\enda\eeq
We here used the fact that $\langle \bfn_t,\bfn\rangle =0$, because $\bfn$ is a unit vector.
In view of (\ref{xn}), we conclude
\bel{xe7}
\beta^\ve(t,\xi)- \beta(t,\xi)~=~
\left\{ \bega{cl} \ov\beta -\beta(\tau,\xi_0)+ \O(1)\cdot\ve
\qquad &\hbox{if}\quad \xi\in S_1\,,\\[1mm]
\O(1) \qquad &\hbox{if}\quad\xi\in S_2\,,\\[1mm]
0 \qquad &\hbox{if}\quad \xi\in S_3\,.\enda \right.
\eeq
For $t\in [\tau-\ve^4, \tau]$, by (\ref{xe6}) and (\ref{xe7}) we obtain
\bel{Et1}\bega{rl}
\E^\ve(t) - \E(t)&\ds=~\int_S \Big[ E\bigl( \beta^\ve(t,\xi)\bigr) \,\bigl| x^\ve_\xi(t,\xi)\bigr|
- E\bigl( \beta(t,\xi)\bigr) \,\bigl| x_\xi(t,\xi)\bigr|\Big]\, d\xi\\[4mm]
&\ds =~\int_{|\xi-\xi_0|<\ve} \bigl[E(\ov\beta) - E(\beta(\tau, \xi_0)\bigr]\, \bigl| x_\xi(\tau,\xi_0)\bigr|
\,d\xi + \O(1)\cdot \ve^2\\[4mm]
&=~2\ve \,  \bigl[E(\ov\beta) - E(\beta(\tau, \xi_0)\bigr]\, \bigl| x_\xi(\tau,\xi_0)\bigr|
+ \O(1)\cdot \ve^2.\enda\eeq
Integrating over time, we finally obtain
\bel{eee}
\int_{\tau-\ve^4}^\tau \Big[\phi \bigl(\E^\ve(t)\bigr) -\phi\bigl(\E(t)\bigr)\Big]\, dt~=~
\phi'\bigl( \E( \tau)\bigr)\cdot 2\ve^5
\bigl[E(\ov\beta) - E(\beta(\tau, \xi_0)\bigr]\, \bigl| x_\xi(\tau,\xi_0)\bigr| + \O(1)\cdot \ve^6.
\eeq
\v
{\bf 4.} In this step we compute the difference $\E^\ve(t)-\E(t)$ in the control effort, 
during the time interval
$t\in [\tau, T]$.
For $(t,\xi) \in [\tau, T]\times S$,  the solution $X$ of (\ref{Xt2})-(\ref{Xt3}) will satisfy different bounds
over the above three sets:
\bel{Xxi}
X_\xi(t,\xi)~=~
\left\{\bega{cl}
\quad \O(1) \quad&\hbox{if}\quad \xi \in S_1\,,\\
\O(1)\cdot \ve^{-2} \quad&\hbox{if}\quad \xi \in S_2\,,\\
\quad 0  \quad&\hbox{if}\quad  \xi \in S_3\,.
\enda\right.
\eeq
Therefore
\bel{dif1}
x^\ve_\xi(t,\xi)- x_\xi(t,\xi)~=~ \begin{cases}
 \O(1)\cdot  \ve^4\quad&\hbox{if}\quad  \xi \in S_1\,,\\
\O(1) \cdot \ve^2\quad&\hbox{if}\quad  \xi \in S_2\,,\\
\quad 0 \quad  \quad&\hbox{if}\quad  \xi \in S_3\,.
\end{cases}
\end{equation}	
The change in the normal speed is computed by
\bel{aa1}\bega{rl}
\beta^\ve -\beta&=~\la x^\ve_t, \,\bfn^\ve\ra - \la x_t, \,\bfn\ra
\\[3mm]
&=~ \ve^4 X_t + \Big\langle \beta\bfn + \ve^4 X_t \bfn + \ve^4 X \bfn_t \,,~\bfn^\ve - \bfn\Big\rangle.
\enda
\eeq
We observe that $\bfn^\ve$ and $\bfn$ are unit vectors.  For $\xi\in S_1$, 
by (\ref{xn}) the angle  between them is
$$\theta_\ve~=~\O(1)\cdot |\bfn^\ve- \bfn|~=~\O(1)\cdot \ve^4.$$
Therefore,
\bel{nn1}\la \bfn, \, \bfn^\ve - \bfn\ra~=~\cos \theta_\ve - 1 ~=~\O(1) \cdot \theta_\ve^2~=~
\O(1)\cdot \ve^8.\eeq
Similarly, when $\xi\in S_2$, by (\ref{xn}) it follows 
\bel{nn2}\la \bfn, \, \bfn^\ve - \bfn\ra~=~
\O(1)\cdot \ve^4.\eeq
Combining (\ref{aa1}) with  (\ref{nn1})-(\ref{nn2}), one obtains the bounds
\bel{aa2}R(t,\xi)~\doteq~
\beta^\ve(t,\xi)-\beta(t,\xi) - \ve^4 X_t(t,\xi)~=~
\left\{\bega{cl}
 \O(1) \cdot \ve^8\quad&\hbox{if}\quad \xi \in S_1\,,\\
\O(1)\cdot \ve^4 \quad&\hbox{if}\quad \xi \in S_2\,,\\
0  \quad&\hbox{if}\quad  \xi \in S_3\,.
\enda\right.
\eeq
We now compute 
\bel{h1}\bega{rl}
  \bigl|x^\ve_\xi\bigr|- |x_\xi|&=\ds~\Big\langle ( x + \ve^4 X \bfn)_\xi\, ,~
( x + \ve^4 X \bfn)_\xi\Big\rangle^{1/2}-|x_\xi|\\[3mm]&\ds=~ \ve^4 \,\left\langle {x_\xi\over |x_\xi|} \,,~X_\xi\bfn + X \bfn_\xi\right\rangle+ \O(1)\cdot \ve^8\\[4mm]
&=~\ds  - \ve^4 \la \bfn^\perp,\bfn_\xi\ra X+\O(1)\cdot \ve^8\\[4mm]
&=~-\ve^4 \omega |x_\xi|\, X+\O(1)\cdot \ve^8\, ,\enda
\eeq
where $\omega$ denotes  the signed curvature, as in (\ref{curv}).

Combining (\ref{h1}) with the evolution equation (\ref{Xt2}) for the perturbation $X$, and using (\ref{aa2}),
we obtain
\bel{Ep1}\bega{l}
E(\beta^\ve) \,|x^\ve_\xi| - E(\beta) \,|x_\xi| \\[4mm]
\qquad =~\Big[ E(\beta^\ve)- E(\beta)\bigr] \,|x_\xi| + E(\beta)\,\big(  |x^\ve_\xi|- |x_\xi|\bigr)
+ \Big[ E(\beta^\ve)- E(\beta)\bigr]\,\big(  |x^\ve_\xi|- |x_\xi|\bigr)\\[4mm]
\qquad =~E'(\beta) [\beta^\ve-\beta] \, |x_\xi| + \O(1)\cdot |\beta^\ve-\beta|^2 +
E(\beta)\,\big(  |x^\ve_\xi|- |x_\xi|\bigr) + \O(1) \cdot |\beta^\ve -\beta|\,\big(  |x^\ve_\xi|- |x_\xi|\bigr)\\[4mm]
\qquad  =~E'(\beta) (\ve^4 X_t+R)  \, |x_\xi| + \O(1)\cdot |\beta^\ve-\beta|^2 +
E(\beta)\,\big( -\ve \omega  |x_\xi| X +\O(1)\cdot \ve^8\bigr) \\[4mm]
\qquad\qquad + \O(1) \cdot |\beta^\ve -\beta|\,\big(  |x^\ve_\xi|- |x_\xi|\bigr),
\enda
\eeq
\bel{Ep2}
E\bigl(\beta^\ve(t,\xi)\bigr) \, \bigl|x^\ve_\xi(t,\xi)\bigr| - E\bigl(\beta(t,\xi)\bigr) \bigl|x_\xi(t,\xi)\bigr|
~=~\left\{\bega{cl}
 \O(1) \cdot \ve^8\quad&\hbox{if}\quad \xi \in S_1\,,\\
\O(1)\cdot \ve^4 \quad&\hbox{if}\quad \xi \in S_2\,,\\
0  \quad&\hbox{if}\quad  \xi \in S_3\,.
\enda\right.
\eeq
Integrating over the whole set $S$, for every $t\in [\tau,T]$ we thus obtain
\bel{Ep3}
\E^\ve(t)-\E(t)~=~\O(1)\cdot \ve^8\, \meas(S_1) + \O(1)\cdot \ve^4\,\meas(S_2)
~=~\O(1)\cdot \ve^6,\eeq
and finally
\bel{ee3}
\int_\tau^T \bigl[ \E^\ve(t) -\E(t)\bigr]\, dt~=~ \O(1)\cdot \ve^6.
\eeq
In other words, by  the identity
(\ref{Xt2}) and  the choice of the  function $x^\ve$ in (\ref{xe2}), 
the change in the cost of the control $\beta$ over the remaining time interval 
$[\tau,T]$ vanishes, to higher order.
\v
{\bf 5.} It remains to estimate the change in the running cost and in the terminal cost
for the perturbed strategies.   We compute
\bel{run-C}
\bega{l}
\ds\left(\int_{\tau-\ve^4}^\tau+ \int_\tau^T\right)\Big[ m_2\bigl(\Omega^\ve (t)\bigr)- m_2\bigl(\Omega(t)\bigr)\Big] \, dt\\[4mm]
\ds\qquad =~\O(1)\cdot\ve^9 - \int_\tau^T  \int_{S_1\cup S_2}\Big(  \ve^4X(t,\xi) + \O(1)\cdot \ve^5
\Big)
 \,\bigl| x_\xi(t,\xi)\bigr|\, d\xi \, dt  \\[4mm]
\ds \qquad =~- \ve^4 
\int_\tau^T  \int_{S_1}  X(t,\xi)
 \,\bigl| x_\xi(t,\xi)\bigr|\, d\xi \, dt  + \O(1)\cdot \ve^6.
\enda
\eeq
Moreover, the change in the final area can be estimated as
\bel{fin-C}\bega{l}
\ds m_2\bigl(\Omega^\ve (T)\bigr)- m_2\bigl(\Omega(T)\bigr)~=~ - \int_{S_1\cup S_2}\Big(  \ve^4X(T,\xi) + \O(1)\cdot \ve^5
\Big)
 \,\bigl| x_\xi(T,\xi)\bigr|\, d\xi\\[4mm]
 \qquad\ds =~- \ve^4 
\int_\tau^T  \int_{S_1}  X(T,\xi)
 \,\bigl| x_\xi(T,\xi)\bigr|\, d\xi \, dt  + \O(1)\cdot \ve^6.
\enda
\end{equation}
We now claim that, if  the adjoint variable $Y$ satisfies (\ref{Ydt})-(\ref{YT}), then 
the change in cost can be computed by
\bel{Y0}
c_1 \int_\tau^T \int_{S_1} \bigl| x_\xi(t,\xi)\bigr| X(t,\xi)\, d\xi d t +c_2  \int_{S_1} \bigl| x_\xi(T,\xi)\bigr| X(T,\xi)\, d\xi~=~\int_{S_1} \bigl| x_\xi(\tau,\xi)\bigr| X(\tau,\xi)\, Y(\tau,\xi)\, d\xi.\eeq
Indeed, this is trivially true  when $\tau = T$, because in this case $Y(T,\xi)=c_2$.
Moreover, differentiating (\ref{Y0}) w.r.t.~time $\tau$,
by  (\ref{Xt2}) and (\ref{Ydt})  we obtain
\bel{Y3}\bega{l}
\ds {d\over d\tau} 
\int_{S_1} \bigl| x_\xi(\tau,\xi)\bigr| \,X(\tau,\xi)\, Y(\tau,\xi)\, d\xi\\[4mm]
\qquad \ds=~\int_{S_1} \left( - \omega(\tau,\xi)\, \beta(\tau,\xi) +     {E\bigl(\beta(\tau,\xi)\bigr)\over E'\bigl(\beta(\tau,\xi)\bigr)}\cdot\omega(\tau,\xi)   \right) \bigl|x_\xi(\tau,\xi)\bigr|   X(\tau,\xi)\, Y(\tau,\xi)\, 
d\xi \\[4mm]
\qquad\qquad \ds +\int_{S_1} \bigl| x_\xi(\tau,\xi)\bigr| X(\tau,\xi)\, Y_\tau(\tau,\xi)\, d\xi \\[4mm]
\qquad \ds =~ -c_1 \int_{S_1} \bigl| x_\xi(\tau,\xi)\bigr| X(\tau,\xi)\, d\xi ~=~{d\over d\tau}
 \int_\tau^T \int_{S_1}  c_1\bigl| x_\xi(t,\xi)\bigr| X(t,\xi)\, d\xi d t\,
.\enda
\eeq
This shows that the identity (\ref{Y0}) holds for every $\tau$.

Together with (\ref{run-C})-(\ref{fin-C}),  from (\ref{Y0}) we obtain
\bel{Y4} 
\bega{l}
\ds c_1\,\left(\int_{\tau-\ve^4}^\tau+ \int_\tau^T\right)\Big[ m_2\bigl(\Omega^\ve (t)\bigr)- m_2\bigl(\Omega(t)\bigr)\Big] \, dt+ c_2 \, m_2\bigl(\Omega^\ve (T)\bigr)- m_2\bigl(\Omega(T)\bigr)   \\[4mm]
\qquad =~\ds
-\ve^4 \,\int_{S_1} \bigl| x_\xi(\tau,\xi)\bigr| \,X(\tau,\xi)\, Y(\tau,\xi)\, d\xi+
 \O(1)\cdot \ve^6\\[4mm]
\qquad =~-2\ve^5 \bigl| x_\xi(\tau,\xi_0)\bigr| \,X(\tau,\xi_0)\, Y(\tau,\xi_0) + \O(1)\cdot \ve^6
\\[4mm]
\qquad =~-2\ve^5 \bigl| x_\xi(\tau,\xi_0)\bigr| \,\bigl( \ov\beta - \beta(\tau, \xi_0)\bigr) Y(\tau,\xi_0) + \O(1)\cdot \ve^6

\enda
\eeq
\v
{\bf 6.} By assumption, the cost of the
perturbation cannot be lower than the original cost.
In view of (\ref{eee}), (\ref{ee3}), and (\ref{Y4}), this implies
\bel{comp}\phi'\bigl(\E(\tau)\bigr)\cdot
2\ve^5 
\bigl[E(\ov\beta) - E(\beta(\tau, \xi_0)\bigr]\, \bigl| x_\xi(\tau,\xi_0)\bigr| 
- 2\ve^5 \bigl| x_\xi(\tau,\xi_0)\bigr| \, \,\bigl( \ov\beta - \beta(\tau, \xi_0)\bigr)  Y(\tau,\xi_0) + \O(1)\cdot \ve^6~\geq~0\eeq
for every $(\tau,\xi)\in \,]0,T[\,\times S$ and every speed $\ov\beta\geq \beta_0$.
Since $\ve>0$ can be taken arbitrarily small, from (\ref{comp}) we deduce
\bel{E33}
\phi'\bigl(\E(\tau)\bigr)\bigl[E(\ov\beta) - E(\beta(\tau, \xi_0)\bigr] -
\bigl( \ov\beta - \beta(\tau, \xi_0)\bigr)  Y(\tau,\xi_0)~\geq~0\eeq
for every $\ov\beta\geq \beta_0$.   This proves (\ref{max1}).

Finally, by continuity the same conclusion remains valid also for $t=0$ or $t=T$. 
\endproof

\begin{remark}\label{r:61} {\rm 
The adjoint variable $Y>0$ introduced at 
(\ref{Ydt})-(\ref{YT}) can be interpreted as a ``shadow price".
Namely (see Fig.~\ref{f:sc21}),  assume that  at time $\tau$ an external contractor offered to remove the contamination
from a neighborhood of the point $x(\tau, \xi_0)$, thus replacing the set $\Omega(\tau)$ with a smaller set $\Omega^\ve(\tau)$, 
 at a price of $Y(\tau,\xi_0)$ per unit area.   
In this case, accepting or refusing the offer would make no difference in the total cost.
 }
\end{remark} 

\begin{remark}\label{r:62} {\rm 
In order to derive the necessary conditions  (\ref{max1}),
we assumed that the parameterization $(t,\xi)\mapsto x(t,\xi)$
had $\C^2$ regularity.   In several applications, this map is continuously differentiable, but
only piecewise $\C^2$.  In particular (see the example in Section~\ref{s:8}), the curvature
$\omega(t,\xi)$ may only be piecewise continuous.  It is worth noting that the proof of 
Theorem~\ref{t:51} remains valid also in this slightly more general setting.
 }
\end{remark}

\section{The case with constraint on the total effort}
\label{sec:6}
\setcounter{equation}{0}
We now consider again the optimization problem 
 {\bf (OP)}, but in the case where the cost function $\phi$ is given by (\ref{EP}).
 This is equivalent to an optimization problem with constraint on the total effort:
\bel{cost6}
\hbox{minimize:}\qquad \J(\Omega)~=~
c_1\int_0^T m_2(\Omega(t))\, dt + c_2\,\meas\bigl(\Omega(T)\bigr),\eeq
subject to
\bel{ETM}
\E(t)~\doteq~ \int_{\partial\Omega(t)}E(\beta(t,x))\, d\sigma~\leq~M\qquad
\hbox{
for every}~~~t\in [0,T].\eeq
This leads to a somewhat different set of necessary conditions.

\begin{theorem} \label{t:61}  Let $E$ satisfy the assumptions {\bf (A1)}.
Assume that $t\mapsto \Omega(t)$ provides an optimal solution to     
(\ref{cost6})-(\ref{ETM}).   
Let $\xi\mapsto x(t,\xi)$ be a $\C^2$ parameterization of the boundary
of the set $\Omega(t)$, satisfying the regularity properties {\bf (A3)}.
Call $Y=Y(t,\xi)$ the adjoint function constructed at (\ref{Ydt})-(\ref{YT}).

Then, for every 
$t\in [0,T]$ one has $\E(t)=M$.
Moreover, there exists a scalar function $t\mapsto \lambda(t)> 0$
such that  the normal velocity $\beta=\beta(t,\xi)$ satisfies
\bel{max16}
\lambda(t) E\big(\beta(t,\xi)\bigr) -Y(t,\xi) \beta(t,\xi)  ~= ~\min_{\beta\geq \beta_0}~
\Big\{\lambda(t) E(\beta)-Y(t,\xi) \beta  \Big\}.\eeq
\end{theorem}

{\bf Proof.} {\bf 1.} To prove the first statement, we argue by contradiction.
If $\E(\tau)<M$, by continuity we can assume
$\E(t)<M$ for all $t$ in a neighborhood of $\tau$.
Then we can choose any $\xi_0\in S$ and any constant $\ov\beta>\beta(\tau, \xi_0)$.
For $\ve>0$ small, we define the perturbed strategy $x^\ve$ on $[\tau-\ve^4, \tau]$ 
as in (\ref{Xt3}), (\ref{xe3}).
Since we are changing the inward velocity $\beta$ only when  $|\xi-\xi_0|<\ve+\ve^2$,
for $\ve>0$ small enough the corresponding total effort will satisfy
$$\E^\ve(t)~\leq~\E(t) + C\,\ve~<~M \qquad\forall t\in [\tau-\ve^4, \tau].$$
Notice that this perturbation satisfies
$$\Omega^\ve(t) \subseteq\Omega(t)\qquad \forall t\in [0, \tau].$$
Moreover, the inclusion is strict for $\tau-\ve<t\leq\tau$.

For $t\in [\tau,T]$ we now define
$$\Omega^\ve(t)~=~\Omega(t)\cap B\Big( \Omega^\ve(\tau), |\beta_0|(t-\tau)\Big).$$
This yields
\bel{E<M}\E^\ve(t)~\leq~\E(t) ~\leq~M \qquad\forall t\in [\tau, T].\eeq

Indeed, at a.e.~boundary point $x\in \partial\Omega^\ve(t)$, two cases can arise:
\begi
\item[(i)] $x\in \partial\Omega^\ve(t)\cap \partial \Omega(t)$. Then the inward normal speed 
at $x$
is the same as in the original solution.
\item[(ii)] $x\in \partial\Omega^\ve(t)\setminus \partial \Omega(t)$, so that
$$d\bigl( x\,;~\Omega^\ve(\tau)\bigr) ~=~ |\beta_0|(t-\tau).$$
In this case, the inward normal speed is precisely $\beta_0$, and this comes at zero cost.
\endi
Combining the two above cases, we obtain (\ref{E<M}).

In conclusion, we obtained  an admissible motion $t\mapsto\Omega^\ve(t)\subseteq\Omega(t)$, which achieves the strict inequality
$$\int_{\tau-\ve^4}^\tau \Omega^\ve(t)\, dt ~<~\int_{\tau-\ve^4}^\tau \Omega(t)\, dt .$$
This contradicts the optimality of $\Omega(\cdot)$.
\v
{\bf 2.} To prove the second statement, we need to find $\lambda(t)>0$ for which
the (\ref{max16}) holds.

Fix a time $\tau$, and consider any two points $\xi_1,\xi_2\in S$ where the control is
active:
$$\beta(\tau, \xi_1)\,>\,\beta_0,\qquad \beta(\tau, \xi_2)\,>\,\beta_0.$$
We claim that for $i=1,2$ the ratios 
${Y(\tau, \xi_i) / E'\bigl(\beta(\tau,\xi_i)\bigr)}$
must be equal.    
If they are not, 
assuming that 
\bel{rat} 
{Y(\tau, \xi_1) \over E'\bigl(\beta(\tau,\xi_1)\bigr)}~>~  
 {Y(\tau, \xi_2) \over E'\bigl(\beta(\tau,\xi_2)\bigr)}~>~0\,,\eeq
we will obtain a contradiction.  

Indeed, recalling that  $Y>0$, by (\ref{rat}) we deduce the existence
of $\delta_1,\delta_2>0$ small enough such that
\bel{YY0} \delta_1 E'(\tau, \xi_1) ~>~ \delta_2 E'(\tau, \xi_2) ,\eeq
\bel{Yd}\delta_1 Y(\tau,\xi_1)~>~
\delta_2 Y(\tau,\xi_2)~>~0.\eeq
Since the function $E$ is continuously differentiable, by possibly shrinking the values of $\delta_1, \delta_2$ 
while keeping the ratio $\delta_1/\delta_2$ constant,
by (\ref{YY0}) we obtain
\bel{YY} \bega{l}\ds \bigl|x_\xi(\tau,\xi_1)\bigr| \left[E\left(\beta(\tau, \xi_1) 
+{\delta_1\over \bigl|x_\xi(\tau,\xi_1)\bigr| }\right) - 
E\bigl(\beta(\tau, \xi_1) \bigr) \right]\\[4mm]
 \qquad\ds <~ \bigl|x_\xi(\tau,\xi_2)\bigr| \left[ E
\bigl(\beta(\tau, \xi_2) \bigr) -
E\left(\beta(\tau, \xi_2) -{\delta_2\over \bigl|x_\xi(\tau,\xi_1)\bigr| }\right)\right] ,\enda \eeq
For $\ve>0$ small 
we construct a  perturbation $x^\ve(t,\xi)$ as in the proof of Theorem~\ref{t:51},
but taking place simultaneously over the two disjoint intervals 
$$I_1\cup I_2~=~\{ \xi\,;~|\xi-\xi_1| < \ve+\ve^2\} \cup \{ \xi\,;~|\xi-\xi_2| < \ve+\ve^2\}.$$
At time $\tau$ we define
\bel{Xt6} X^\ve(\tau, \xi)~=~{\delta_1\over \bigl| x_\xi(\tau,\xi_1)\bigr|}  \, \vp_\ve(\xi-\xi_1) 
- {\delta_2\over \bigl| x_\xi(\tau,\xi_2)\bigr|}\,\vp_\ve(\xi-\xi_2) .\eeq
On the remaining interval $[\tau, T]$, we define $X^\ve$ to be the solution of 
\bel{Xt22}
X^\ve_t(t,\xi)~=~{E\bigl(\beta(t,\xi)\bigr)\over 
E'\bigl(\beta(t,\xi)\bigr)}\cdot \omega(t,\xi)\, X^\ve(t,\xi) - \ve^{2/3}\,\bigl|X^\ve(t,\xi)\bigr|, \eeq
with initial data (\ref{Xt6})  at $t=\tau$.
Notice that, compared with (\ref{Xt2})-(\ref{Xt3}), here the construction of $X^\ve$ includes a further $\ve$-perturbation. 
This is needed, in order to guarantee that the total effort remains $\leq M$ at all times.

Similarly to the proof of Theorem~\ref{t:51}, for all $\xi\in S$ we now define  \bel{xe33}
x^\ve(t,\xi)~=~\left\{
\bega{cl} x(t,\xi) &\hbox{if}~~t\in [0, \tau-\ve^4],\\[3mm]
x(t,\xi)+ \bigl[ \ve^4 - (\tau-t)\bigr]\,X^\ve(\tau,\xi)\bfn(t,\xi)\qquad &\hbox{if}\quad t\in [\tau-\ve^4, \,\tau],
\\[3mm]
x(t,\xi)+ \ve^4 X^\ve(\tau,\xi)\bfn(t,\xi)&\hbox{if}~~t\in [\tau,T].\enda\right.
\eeq
\v
{\bf 3.} In this step we prove that, for $ \ve >0 $ sufficiently small
and every  $t\in \,]\tau-\ve^4, \tau[\,$, the instantaneous effort
satisfies 
\bel{ieff}\E^\ve(t)~\leq~ \E(t)~\leq~ M.\eeq   For notational convenience, define
%
%
%
%
\begin{equation}
S_1^i \,\doteq\,\{ \xi\in S\,;~|\xi-\xi_i| \le \ve\}, \qquad S_2^i\,\doteq\, \{ \xi\in S \,;~\ve<|\xi-\xi_i| < \ve+\ve^2\},\qquad
i=1,2
\end{equation}
As in (\ref{xe6}), one has
\begin{equation}
x^\ve_\xi(t,\xi)- x_\xi(t,\xi)~=~ \begin{cases}
\O(1)\cdot  \ve^4\quad&\hbox{if}\quad  \xi \in S_1^i\,,\\
\O(1) \cdot \ve^2\quad&\hbox{if}\quad  \xi \in S_2^i\,,\\
\quad 0 \quad  \quad&\hbox{otherwise.}
\end{cases}
\end{equation}
As in (\ref{xn}), the same bounds  hold for  $\bfn^\ve - \bfn$. 

For $\tau-\ve^4<t< \tau$, the same arguments used at (\ref{xe7})
show that
the inward normal velocity satisfies
\bel{614}
\beta^\ve(t,\xi)- \beta(t,\xi)~=~
\left\{ \bega{cl} \ds {\delta_1\over \bigl| x_\xi(\tau,\xi_1)\bigr|}+ \O(1)\cdot\ve
\qquad &\hbox{if}\quad \xi\in S_1^1\,,\\[4mm]
 \ds -{\delta_2\over \bigl| x_\xi(\tau,\xi_2)\bigr|} + \O(1)\cdot\ve
\qquad &\hbox{if}\quad \xi\in S_1^2\,,\\[4mm]
\O(1) 
\qquad &\hbox{if}\quad\xi\in S_2^1\cup S_2^2\,,\\[2mm]
0\qquad &\hbox{otherwise.}\enda \right.
\eeq
For 
$\tau-\ve^4<t<\tau$, using (\ref{614}) and (\ref{YY}), we now compute
\bel{per2}\bega{l}
\E^\ve(t) - \E(t)~\ds=~\int_S \Big[ E\bigl( \beta^\ve(t,\xi)\bigr) \,\bigl| x^\ve_\xi(t,\xi)\bigr|
- E\bigl( \beta(t,\xi)\bigr) \,\bigl| x_\xi(t,\xi)\bigr|\Big]\, d\xi\\[4mm]
\qquad \ds =~\ve\Bigg\{\bigl|x_\xi(\tau,\xi_1)\bigr| \left[E\Big(\beta(\tau, \xi_1) +{\delta_1\over \bigl|x_\xi(\tau,\xi_1)\bigr| }\Big) - 
E\bigl(\beta(\tau, \xi_1) \bigr) \right]\\[4mm]
\qquad\qquad\qquad \ds -
\bigl|x_\xi(\tau,\xi_2)\bigr| \left[ E
\bigl(\beta(\tau, \xi_2) \bigr) -
E\Big(\beta(\tau, \xi_2) -{\delta_2\over \bigl|x_\xi(\tau,\xi_1)\bigr| }\Big)\right] \Bigg\}\\[4mm]
\qquad \quad \ds+ \int_{S_1^1\cup S_1^2} \O(1)\cdot\ve\, d\xi + \int_{S_2^1\cup S_2^2} \O(1)\, d\xi 
\\[4mm]
\qquad < ~0,\enda
\eeq
provided that $\ve>0$ is sufficiently small.
Indeed, the first  term on the right hand side of (\ref{per2}) is strictly negative, while the last two terms are of order $\O(1)\cdot \ve^2$. 
\v
{\bf 4.}
Next, we estimate the change in the control effort $\E(t)$ for $\tau<t<T$. Here the computations are very similar
to the ones in step {\bf 4} of the proof of Theorem~\ref{t:51}.  Because of the additional term on the right hand side of 
(\ref{Xt22}), the bounds
(\ref{aa2}) are now replaced by 
\bel{aa22}R(t,\xi)~\doteq~
\beta^\ve(t,\xi)-\beta(t,\xi) - \ve^4 X_t(t,\xi)~=~
\left\{\bega{cl}
 \O(1) \cdot \ve^8 - \ve^4 \ve^{2/3} \bigl|X(t,\xi)\bigr|\quad&\hbox{if}\quad \xi \in S_1^1\cup S_1^2\,,\\
\O(1)\cdot \ve^4 \quad&\hbox{if}\quad \xi \in S_2^1\cup S_2^2\,,\\
0  \quad&\hbox{otherwise.}
\enda\right.
\eeq
In turn, the estimates (\ref{Ep2}) are replaced by
\bel{Ep22}\bega{l}
E\bigl(\beta^\ve(t,\xi)\bigr) \, \bigl|x^\ve_\xi(t,\xi)\bigr| - E\bigl(\beta(t,\xi)\bigr) \bigl|x_\xi(t,\xi)\bigr|\\[4mm]
\ds\qquad =~\left\{\bega{cl}  - \ve^4 \ve^{2/3} \bigl|X(t,\xi)\bigr| E'\bigl(\beta(t,\xi)\bigr)+
 \O(1) \cdot \ve^8\quad&\hbox{if}\quad \xi \in S_1^1\cup S_1^2\,,\\
\O(1)\cdot \ve^4 \quad&\hbox{if}\quad \xi \in S_2^1\cup S_2^2\,,\\
0  \quad&\hbox{otherwise.}
\enda\right.\enda
\eeq
We now integrate over the whole set $S$.   Observing that, for $\ve>0$ small, the 
function $|X(t,\xi)\bigr|$ remains uniformly positive over the set $ (S_1^1\cup S_1^2)\times [\tau,T]$, 
for a suitable constant $c_0>0$
and all $\tau<t<T$ we obtain
\bel{Ep33}\bega{rl}
\E^\ve(t)-\E(t)&\ds=~- \int_{S_1^1\cup S_1^2}  \ve^4 \ve^{2/3} \bigl|X(t,\xi)\bigr| E'\bigl(\beta(t,\xi)\bigr) \, d\xi \\[4mm]
&\qquad 
+   \O(1)\cdot \ve^8\, \meas(S_1^1\cup S_1^2) + \O(1)\cdot \ve^4\,\meas(S_2^1\cup S_2^2)\\[3mm]
&\leq - c_0 \ve^5 \ve^{2/3} + \O(1)\cdot \ve^6~<~0.\enda\eeq
\v
{\bf 5.} It remains to estimate the change in the running cost and in the terminal cost
for the perturbed strategies.  We notice that the formulas 
(\ref{run-C})-(\ref{fin-C}) remain valid.

We now consider the auxiliary function 
$\Tilde X: [\tau,T]\times S\mapsto \R$, defined to be the solution to the Cauchy problem
\bel{TXX}
\Tilde X_t(t,\xi)~=~{E\bigl(\beta(t,\xi)\bigr)\over 
E'\bigl(\beta(t,\xi)\bigr)}\cdot \omega(t,\xi)\, \Tilde X(t,\xi), \eeq
\bel{TX0}\Tilde X(\tau, \xi)~=~{\delta_1\over \bigl| x_\xi(\tau,\xi_1)\bigr|}  \, \vp_\ve(\xi-\xi_1) 
- {\delta_2\over \bigl| x_\xi(\tau,\xi_2)\bigr|}\,\vp_\ve(\xi-\xi_2) .\eeq
A comparison with (\ref{Xt6})-(\ref{Xt22}) yields
\bel{XTX} \bigl| X^\ve(t,\xi)-\Tilde X(t,\xi)\bigr|~\leq~C\, \ve^{2/3}\eeq
for all $t,\xi$.
We now call $\Omega^\ve(t)$ and $\Tilde\Omega^\ve(t)$ respectively the sets 
corresponding to the perturbations 
$$x^\ve(t,\xi)~=~x(t,\xi) + \ve^4 X^\ve(t,\xi) \bfn(t,\xi), \qquad
\tilde x^\ve(t,\xi)~=~x(t,\xi) + \ve^4 \Tilde X(t,\xi) \bfn(t,\xi).$$
In view of (\ref{Y4}), with $X$ replaced by $\Tilde X$, recalling (\ref{Yd}) we conclude
\bel{Y44} 
\bega{l}
\ds c_1\,\left(\int_{\tau-\ve^4}^\tau+ \int_\tau^T\right)\Big[ m_2\bigl(\Omega^\ve (t)\bigr)- m_2\bigl(\Omega(t)\bigr)\Big] \, dt+ c_2 \, \Big[m_2\bigl(\Omega^\ve (T)\bigr)- m_2\bigl(\Omega(T)\bigr) \Big]  \\[4mm]
\qquad =~\ds  \O(1)\cdot \ve^4\, \ve\, \ve^{2/3} + c_1 \int_\tau^T\Big[ m_2\bigl(\Tilde 
\Omega^\ve (t)\bigr)- m_2\bigl(\Omega(t)\bigr)\Big] \, dt+ c_2 \, \Big[m_2\bigl(\Tilde\Omega^\ve (T)\bigr)- m_2\bigl(\Tilde \Omega(T)\bigr)\Big] \\[4mm]
\qquad =~\ds
\ds  \O(1)\cdot \ve^{17/3} -
\ve^4 \,\int_{S_1^1\cup S_1^2} \bigl| x_\xi(\tau,\xi)\bigr| \,\Tilde X(\tau,\xi)\, Y(\tau,\xi)\, d\xi+
 \O(1)\cdot \ve^6\\[4mm]
\qquad =~-2\ve^5 \bigg[ \bigl| x_\xi(\tau,\xi_1)\bigr| \,\Tilde X(\tau,\xi_1)\, Y(\tau,\xi_1) + \bigl| x_\xi(\tau,\xi_2)\bigr| \,
\Tilde X(\tau,\xi_2)\, Y(\tau,\xi_2)\bigg] + \O(1)\cdot \ve^{17/3}
\\[4mm]
\qquad =~-2\ve^5\Big( \delta_1 Y(\tau,\xi_1) - \delta_2 Y(\tau,\xi_2)\Big) + \O(1)\cdot \ve^{17/5},
\enda
\eeq
for all $\ve>0$ small enough.
\v
{\bf 6.} The previous steps have established the existence of a function $t\mapsto \lambda(t)>0$ such that
\bel{ratio}E'\bigl(\beta(t,\xi)\bigr)~=~ {Y(t,\xi) \over  \lambda(t)}\, \eeq
at all points where $\beta(t, \xi)>\beta_0$.
Since the function $E(\cdot)$ is convex, we conclude
$$\beta(t,\xi)~=~\argmin_{\beta\geq \beta_0} \Big\{E(\beta) -  {Y(t,\xi) \over \lambda(t)}\beta\Big\}~=~\argmin_{\beta\geq \beta_0} \Big\{\lambda(t) E(\beta) -  Y(t,\xi) \beta\Big\}.$$
This yields (\ref{max16}).\endproof

\section{Optimality conditions at junctions}
\label{s:7}
\setcounter{equation}{0}

In Theorem~\ref{t:41}, the existence of optimal solutions  was proved within a class of 
functions with BV regularity.  On the other hand,  the necessary conditions for optimality 
derived in Theorem~\ref{t:51} require that the sets $\Omega(t)$ have $\C^2$  boundary.
Aim of this section is to partially fill this regularity gap, ruling out certain configurations where  the sets $\Omega(t)$ have 
corners.  Toward this goal, we need to strengthen the assumption {\bf(A1)}, replacing (\ref{Eass}) with the strict inequality
\bel{Ess}  E(\beta) -\beta E'(\beta)~>~0\qquad\qquad\forall \beta>0.\eeq
As explained in Remark~\ref{r:14}, if (\ref{Eass}) holds then a wiggly boundary as shown in Fig.~\ref{f:sc54}
cannot achieve a lower cost, compared with a flat boundary. By imposing the stronger assumption (\ref{Ess}), we
make sure that the wiggly boundary yields a strictly larger cost.
A useful consequence of this assumption is
\begin{lemma}\label{l:71}
Assume that  the effort function $E:\R\mapsto \R_+$ satisfies {\bf (A1)}, with (\ref{Eass}) replaced by (\ref{Ess}). 
Then, for all 
$ \beta>\beta_0$ and 
$\lambda>1$ one has
\bel{Ein} E(\lambda  \beta)~<~\lambda E( \beta).\eeq
\end{lemma}

{\bf Proof.} 
{\bf 1.} Assume first $\beta_0<\beta\leq 0$.  In this case, for $\lambda>1$  we have
$$\lambda \beta \leq \beta\leq 0,\qquad E(\beta)>0,
$$ and hence
$$E(\lambda \beta) ~\leq~E(\beta)~<~\lambda \, E(\beta).$$
\v
{\bf 2.} Next, assume $\beta>0$. If $E(\lambda \beta) \geq \lambda E(\beta)$, a contradiction is obtained as follows.
Choose an intermediate value $\beta_1\in [\beta, \lambda\beta]$ such that 
$$E'(\beta_1)~=~{E(\lambda\beta)-E(\beta)\over (\lambda-1)\beta}\,.$$
By convexity, the graph of $E$ lies below the secant line through the points $\beta, \lambda\beta$. 
Therefore
$$\bega{l} \ds 
E(\beta_1) -\beta_1 E'(\beta_1)~\ds\leq ~E(\beta) + {E(\lambda\beta)-E(\beta)\over (\lambda-1)\beta}\cdot (\beta_1-\beta) 
- \beta_1{E(\lambda\beta)-E(\beta)\over (\lambda-1)\beta}\\[4mm]
\qquad \ds=~E(\beta) - {E(\lambda\beta)-E(\beta)\over (\lambda-1)\beta}\,\beta
~\leq~ E(\beta) - {\lambda E(\beta)-E(\beta)\over (\lambda-1)\beta}\,\beta ~=~0,\enda
$$
reaching a contradiction with (\ref{Ess}).  \endproof

\begin{figure}[ht]
\centerline{\hbox{\includegraphics[width=11cm]{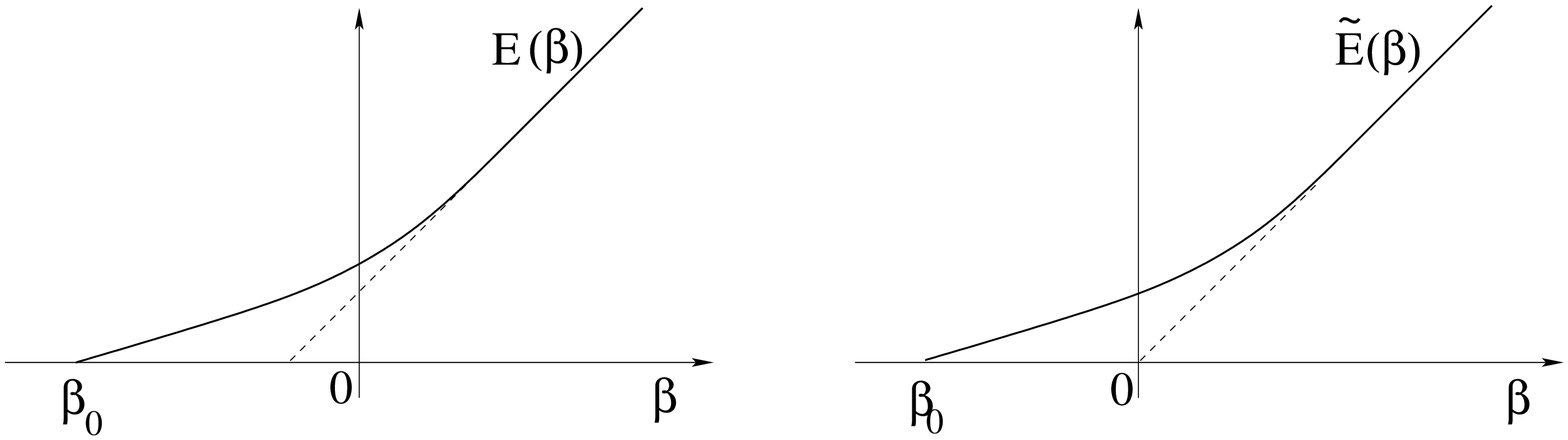}}}
\caption{\small Left: an effort function $E(\beta)$ satisfying the strict inequality (\ref{Ess}). Right: an effort function $\Tilde E(\beta)$
satisfying the assumptions {\bf (A2)} but not (\ref{Ess}).}
\label{f:sc66}
\end{figure}
\begin{figure}[ht]
\centerline{\hbox{\includegraphics[width=13cm]{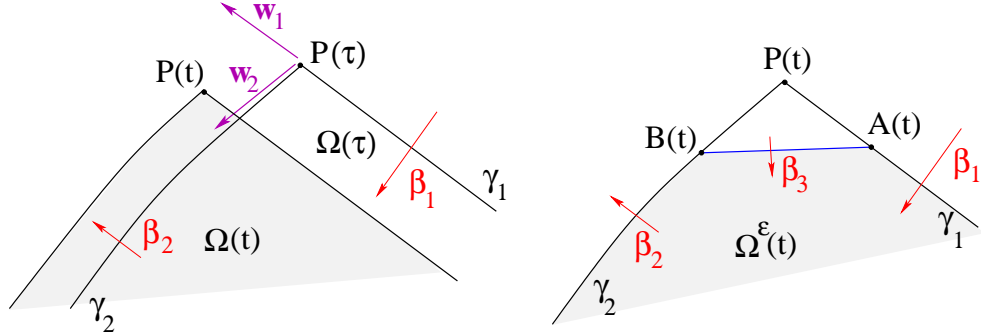}}}
\caption{\small The case of outward corner. Left: the shaded region is the set $\Omega(t)$, for $t\in [\tau, \tau+\delta]$. 
Right: the set $\Omega^\ve$ is obtained from $\Omega(t)$ by removing the triangular
region $\Hat{APB}$.}
\label{f:sc64}
\end{figure}
\begin{figure}[ht]
\centerline{\hbox{\includegraphics[width=13cm]{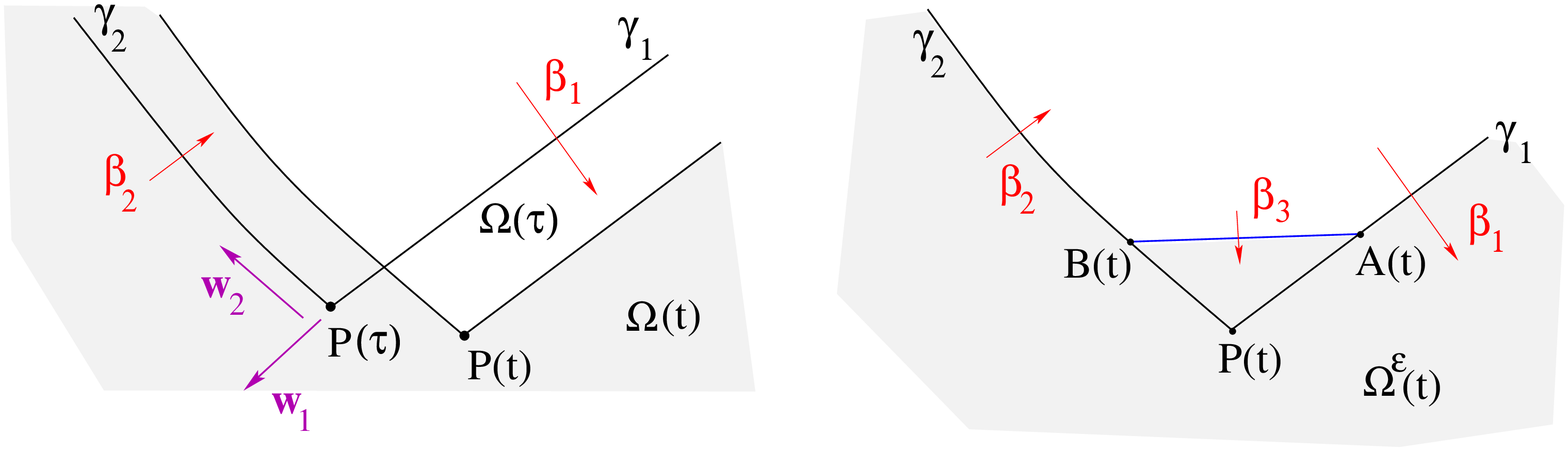}}}
\caption{\small The case of an inward corner, with $\bfw_1\times \bfw_2<0$. 
Left: the shaded region is the set $\Omega(t)$, for $t\in [\tau, \tau+\delta]$. 
Right: the set $\Omega^\ve$ is obtained from $\Omega(t)$ by  adding the triangular
region $\Hat{APB}$.}
\label{f:sc65}
\end{figure}

In the remainder of this section, we shall consider the following situation:
\begi
\item[{\bf (A4)}] {\it
There exists $\tau,\delta_0>0$ such that, for $|t-\tau|<\delta_0$,  the boundary $\partial\Omega(t)$ contains two adjacent arcs $\gamma_1(t,\cdot)$, $\gamma_2(t,\cdot)$ joining at a point $P(t)$ 
at an angle $\theta(t)$.    Each of these arcs admits a $\C^1$  parameterization by arc-length, of the form
\bel{angle}\left\{ \bega{l}s\mapsto \gamma_1(t,s),\qquad s\leq 0,\\[2mm]
s\mapsto \gamma_2(t,s),\qquad s\geq 0,\enda
\right.\qquad\qquad \gamma_1(t,0) = \gamma_2(t,0)
=P(t).\eeq
}\endi
For future reference,
the tangent vectors to the curves $\gamma_1(\tau,\cdot)$ and $\gamma_2(\tau,\cdot)$ at the intersection point
$P(\tau)$ will be denoted by
\bel{bw}\bfw_1~=~\gamma_{1,s}(\tau, 0-)\,,\qquad\qquad \bfw_2~=~\gamma_{2,s}(\tau, 0+)\,.\eeq
Moreover, we call $\bfw^\perp_1,\bfw_2^\perp$ the orthogonal vectors (rotated by $90^o$ counterclockwise).
Notice that the two curves $\gamma_1,\gamma_2$ form an outward corner at $P(\tau)$ if the vector product
satisfies (see Fig.~\ref{f:sc64}, left)
$$\bfw_1\times \bfw_2~\doteq~\la \bfw_1^\perp, \bfw_2\ra~>~0.$$
On the other hand, if $\bfw_1\times \bfw_2<0$ one has an inward corner, as shown in Fig.~\ref{f:sc65}.

As before, we say that the control is {\em active} on a portion of the boundary $\partial \Omega(t)$ if
the inward normal speed is $\beta>\beta_0$.   By {\bf (A1)}, this means that the effort is strictly positive: $E(\beta)>0$.
The main result of this section shows that, for an optimal
motion $t\mapsto \Omega(t)$, non-parallel junctions  cannot be optimal if the control is active on at least
one of the adjacent arcs.  

\begin{figure}[ht]
\centerline{\hbox{\includegraphics[width=11cm]{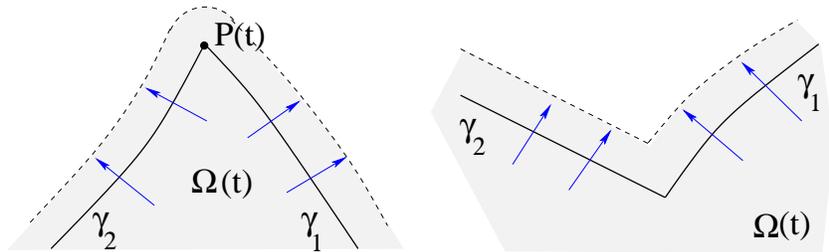}}}
\caption{\small If no control is active, the set $\Omega(t)$ expands with speed $|\beta_0|$ in all directions.
Hence outward corners instantly disappear (left), but inward corners persist (right).}
\label{f:sc60}
\end{figure}

\begin{remark}  {\rm If the control is not active along any of the two arcs $\gamma_1, \gamma_2$,
then the set $\Omega(t)$ expands with speed $|\beta_0|$ all along the boundary, 
in a neighborhood of $P(t)$.  This implies that $\Omega(t)$ satisfies an interior ball condition, hence it can only have
inward corners, as shown in Fig.~\ref{f:sc60}.
}
\end{remark}

\begin{theorem}\label{t:72} Let $E$ satisfy {\bf (A1)}, with (\ref{Eass}) replaced by (\ref{Ess}), and let $\phi:\R_+\mapsto\R_+$ 
be a $\C^1$ function which satisfies {\bf (A2)}. 
Assume that $t\mapsto \Omega(t)$ provides an optimal solution to {\bf (OP)}. 

In the setting described at {\bf (A4)}, if along at least one of the two arcs $\gamma_1, \gamma_2$ the control is active
(i.e., if $\beta>\beta_0$ along the arc), then the two arcs must be tangent at $P(t)$.
\end{theorem}

{\bf Proof.} {\bf 1.} 
Assume that the  two arcs $\gamma_1,\gamma_2$ are not tangent, forming an angle $\theta(t)\neq \pi$ for $|\tau-t|<\delta$. We 
will then construct a perturbed multifunction  $t\mapsto\Omega^\ve (t)$ with a smaller cost.
Given $0<\ve<\!<\delta$, for $t\in [\tau-\ve, \tau+\delta+\ve]$ consider the points
\bel{At} A(t) ~=~\left\{ \bega{cl} \gamma_1\bigl(t, \, \tau-t-\ve\bigr) \quad & \hbox{if} ~~t\in [\tau-\ve,\, \tau],\cr
\gamma_1(t, -\ve) \quad &\hbox{if}~~ t\in [\tau, \,\tau+\delta],\cr
\gamma_1\bigl(t, \, t-\tau -\delta -\ve) \bigr)\quad& \hbox{if}~~ t\in [\tau+\delta,\, \tau+\delta+\ve],\enda\right.\eeq
\bel{Bt}
 B(t) ~=~\left\{ \bega{cl} \gamma_2\bigl(t, \,t-\tau+\ve\bigr) \quad & \hbox{if} ~~t\in [\tau-\ve, \,\tau],\cr
\gamma_2(t, \ve) \quad &\hbox{if}~~ t\in [\tau,\, \tau+\delta],\cr
\gamma_2\bigl(t,\, \tau+\delta+\ve-t \bigr)\quad& \hbox{if}~~ t\in [\tau+\delta,\, \tau+\delta+\ve].\enda\right.\eeq
Observe that $A(t) = B(t) = P(t)$ for $t=\tau-\ve$ and for $t=\tau+\delta+\ve$. 
We now construct a family of perturbed sets $\Omega^\ve(t)$ by the following rules.
\begi
\item[(i)]  For $t\notin [\tau-\ve, \tau+\delta+\ve]$, one has $\Omega^\ve(t)~=~\Omega(t)$.
\item[(ii)] If $\bfw_1\times \bfw_2>0$ (an outward corner), then for $t\in [\tau-\ve, \tau+\delta+\ve]$
the set $\Omega^\ve(t)$ is obtained from $\Omega(t)$ by removing the triangular region with vertices $A(t), P(t), B(t)$,
as shown in Fig.~\ref{f:sc64}.
\item[(iii)] If $\bfw_1\times \bfw_2<0$ (an inward corner), then for $t\in [\tau-\ve, \tau+\delta+\ve]$
the set $\Omega^\ve(t)$ is obtained from $\Omega(t)$ by adding the triangular region with vertices $A(t), P(t), B(t)$,
as shown in Fig.~\ref{f:sc65}.
\endi

\begin{figure}[ht]
\centerline{\hbox{\includegraphics[width=10cm]{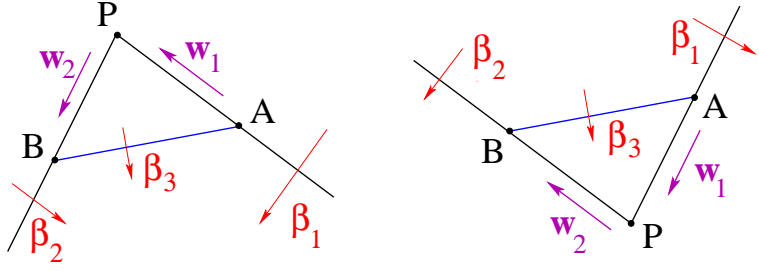}}}
\caption{\small Computing the normal speed $\beta_3$ of the side $AB$, as in (\ref{b123}).}
\label{f:sc47}
\end{figure}

\v
{\bf 2.} To estimate the change in the cost of the new strategy $\Omega^\ve$, the crucial step is
to determine the inward normal speed $\beta_3$ along  the segment $\ov{A(t) B(t)}$. 
Referring to Fig.~\ref{f:sc47} 
consider a triangle with vertices 
$$P(t),\qquad  A(t)\,=\, P(t) - \bfw_1,\qquad  B(t)\,=\, P(t) + \bfw_2\,.$$ 
Call $\beta_1, \beta_2,\beta_3$ respectively the normal speeds of the three sides $AP$, $BP$, and $AB$. 
Knowing the velocity  $\dot P= dP/dt$, these are computed by
\bel{b123}\beta_1~=~\la\dot P, \bfw_1^\perp\ra,\qquad \beta_2~=~\la\dot P, \bfw_2^\perp\ra,
\qquad \beta_3~=~\left\langle \dot P, {(\bfw_1 + \bfw_2)^\perp\over |\bfw_1+\bfw_2|} \right\rangle ~=~{  \beta_1 + \beta_2\over |\bfw_1 +\bfw_2| }.\eeq
Neglecting higher order terms,  during the time interval $[\tau, \tau+\delta]$ the change in the total effort is thus 
computed as
\bel{save}
\bega{rl} \E^\ve(t)-\E(t)& =~\ve\Big( |\bfw_1 +\bfw_2| \, E(\beta_3) -  E(\beta_1) -  E(\beta_2)+\O(1)\cdot (\delta+\ve)\Big)
+o(\ve)\\[3mm]
&=\ds ~\ve \left( |\bfw_1 +\bfw_2| \, E\left( {  \beta_1 + \beta_2\over |\bfw_1 +\bfw_2| }
\right) -  E(\beta_1) -  E(\beta_2)
\right)+\O(1)\cdot \delta\ve+o(\ve).\enda\eeq
We claim that, for $\ve,\delta$ sufficiently small, the right hand side of (\ref{save}) is strictly negative.
Indeed, set 
$$\ov \beta~\doteq~{\beta_1+\beta_2\over 2}~>~\beta_0 \,\qquad\qquad 
\lambda~\doteq ~{2\over |\bfw_1 +\bfw_2|}~>~1 \,.$$ 
Notice that the first inequality follows from the assumption that at least one of the normal speeds $\beta_1,\beta_2$ 
is strictly larger than $\beta_0$.  The second inequality is trivially true because $\bfw_1,\bfw_2$ are non-parallel unit vectors.
By  the strict inequality (\ref{Ein}) and the convexity of $E$ it now follows
\bel{bes}\bega{l}  \ds {1\over 2} \left[{|\bfw_1 +\bfw_2|}  \, E\left( {  \beta_1 + \beta_2\over |\bfw_1 +\bfw_2| }\right) -  E(\beta_1) -  E(\beta_2)\right] =~{1\over \lambda}  E\left( \lambda\cdot { \beta_1 + \beta_2\over 2}\right) - 
{ E(\beta_1) +  E(\beta_2)\over 2} \\[4mm]
\qquad <~\ds E\left( { \beta_1 + \beta_2\over 2}\right) - 
{ E(\beta_1) +  E(\beta_2)\over 2}~\leq~0.
\enda
\eeq
By choosing $0<\ve<\!<\delta$ sufficiently small,  our claim is proved.

\v
{\bf 3.} We now compare the cost of the two strategies $\Omega^\ve(\cdot) $ and $\Omega(\cdot)$, 
for $0<\ve<\!<\delta$ sufficiently small.  For sake of definiteness, we consider an outward corner, 
so that  $\bfw_1\times \bfw_2>0$. The case of an inward corner is entirely similar.
\begi
\item Since $\Omega^\ve(T) =\Omega(T)$, there is no difference in the terminal cost.
\item For every $t\in [\tau-\ve, \tau+\delta+\ve]$ the difference in the area is bounded by
$$\Big|m_2\bigl(\Omega^\ve(t) \bigr)- m_2\bigl(\Omega(t)\bigr)\Big| ~=~\O(1)\cdot \ve^2.$$
Integrating in time, this yields
\bel{cep}
c_1\int_0^T m_2\bigl(\Omega^\ve(t) \bigr)\, dt - c_1\int_0^T m_2\bigl(\Omega(t) \bigr)\, dt~=~\O(1) \cdot \ve^2 (\delta+ 2\ve).\eeq
\item At every time $t\in [\tau-\ve, \,\tau+\delta+\ve]$, the difference in the 
total effort is estimated by an integral of the effort over the segment with endpoints $A(t)$, $B(t)$. Namely
$$\E^\ve(t) - \E(t)~\leq~\int_{\ov{ A(t) \,B(t)}}E\bigl(\beta^\ve(t,x)\bigr)\, d\sigma~ =~\O(1)\cdot |B(t)- A(t)|~=~\O(1)\cdot \ve. $$
Since we are assuming the differentiability of the cost function $\phi$, this implies
\bel{efe}
\left(\int_{ \tau-\ve}^ \tau +\int_{\tau+\delta}^{ \tau+\delta+\ve}\right)  \Big[\phi\bigl(\E^\ve(t)\bigr) - \phi\bigl(\E(t)\bigr)
\Big]\, dt~=~\O(1)
\cdot\ve^2. \eeq

\item  Finally, by the inequalities at (\ref{save})-(\ref{bes}), 
for $t\in [\tau, \tau+\delta]$, the difference in the total effort is bounded by 
$$
\E^\ve(t) - \E(t)~\leq~- \ve \kappa_0 + o(\ve),$$
for some constant $\kappa_0>0$ and all $\delta, \ve>0$ sufficiently small.
Therefore, we can write
\bel{edel}\bega{l} 
\ds \int_\tau^{\tau+\delta}  \Big[\phi\bigl(\E^\ve(t)\bigr) - \phi\bigl(\E(t)\bigr)
\Big]\, dt~\leq ~ \int_\tau^{\tau+\delta} {1\over 2} \phi' \bigl(\E(t)\bigr)\cdot 
\bigl[\E^\ve(t) - \E(t)\bigr ]\,dt \\[4mm]
\ds\qquad\qquad  \leq~ -\int_\tau^{\tau+\delta} {\kappa_0\over 4} \phi' \bigl(\E(t)\bigr)
~\leq~- \kappa_1 \ve\delta\,,
\enda
\eeq
for some constant $\kappa_1>0$.  
\endi
Combining  the estimates (\ref{cep}), (\ref{efe}) and  (\ref{edel}), the difference in the total cost
is estimated as
$$J(\Omega^\ve) - J(\Omega)~\leq~- \kappa_1 \ve\delta + \O(1)\cdot \ve^2~<~0,$$
showing that the original strategy was not optimal.\endproof

\section{Optimal motions determined by the necessary conditions}
\label{s:8}
\setcounter{equation}{0}
In this last section we study the set motions
$t\mapsto \Omega(t)$ that satisfy the necessary conditions derived  earlier.
We focus on the basic case where 
\bel{Ede}
E(\beta)~=~\max\,\{ 0, \, 1+\beta\}.\eeq

In connection with the optimality condition 
\bel{max11}
\lambda(t) E\big(\beta(t,\xi)\bigr) -Y(t,\xi) \beta(t,\xi)  ~= ~\min_{\beta\geq -1}~
\Big\{\lambda(t) E(\beta)-Y(t,\xi) \beta  \Big\},\eeq
three cases need to be considered.

CASE 1: $\lambda(t)-Y(t,\xi)>0$.   In this case $\beta(t,\xi)=-1$. In other words, no effort
is made at the point $x(t,\xi)$.  
Hence the boundary point $x(t,\xi)$ moves in the direction 
of the outer normal with unit speed.

CASE 2: $\lambda(t)-Y(t,\xi)=0$.   In this case, any inward normal speed $\beta(t,\xi)\geq -1$ is compatible with (\ref{max1}).

CASE 3: $\lambda(t)-Y(t,\xi)<0$. This can never happen, because the minimality condition cannot be satisfied.  Formally,
(\ref{max11}) would imply that 
the optimal control is $\beta (t, \xi) = +\infty$. 
%
%
%

Similarly to the case of {\em singular controls}, often encountered in geometric control theory \cite{AS, BoPi, Suss}, 
in Case 2 the pointwise values 
of $\beta(t,\xi)$ are determined not by the minimum principle (\ref{max11}), but 
by the requirement that the function $\xi\mapsto Y(t,\xi)= \lambda(t)$  is independent of $\xi$ 
on the region where the control is active.    
In other words, 
to determine the optimal normal speed $\beta=\beta(t,\xi)$  
we need to use (\ref{Ydt}), and
impose that the right hand side is constant w.r.t.~$\xi$, over the portion of the boundary where the control is active.
When the effort $E$ takes the simple form (\ref{Ede}), the
backward Cauchy problem (\ref{Ydt})-(\ref{YT}) reduces to
\bel{Y6} Y_t(t,\xi)~=~ \omega(t,\xi) Y(t,\xi) - c_1\,, \qquad\qquad Y(T, \xi)=c_2\,.\eeq
In order for $Y(t,\xi) = Y^*(t)$ to be a function of time alone,
this implies the constant curvature condition:
\begi
\item[{\bf (CC)}]
{\it At any time $t\in [0,T]$ the curvature $\omega(t,\cdot)$ must  be constant along the portion of the boundary 
where the control is active.}
\endi
When the cost function $\phi$ is smooth,  this constant value $\lambda(t)$ 
is determined by the scalar
equation (\ref{lamb}).   
On the other hand, when $\phi$ is the function in (\ref{EP}), 
$\lambda(t)$ can be determined by the global constraint
\bel{gc}
\int_S E(\beta(t,\xi))\, \bigl|x_\xi(t,\xi)\bigr|\, d\xi~=~M.\eeq
By the necessary conditions, the control is active at points $x(t,\xi)$ where the 
dual variable $Y(t,\xi)$ has the largest values. These are the points where it is most advantageous to shrink the set $\Omega(t)$.     
By (\ref{Y6}), 
the characteristics $t\mapsto x(t,\xi)$ where the dual function $Y$ grows
faster (going backwards in time) are those where the curvature is maximum, and the control is active.
This leads  to the following 

\begin{conjecture}\label{c:81}
At every time $t\in [0,T]$, the optimal control is  active precisely 
along the portion of the boundary where the curvature is maximum.
\end{conjecture}

The validity of this conjecture will be a topic for future investigation.   Here we conclude with a simple example,
where a motion satisfying the necessary conditions can be computed explicitly.

\begin{example}\label{e:81}
{\rm Assume that the initial set $\Omega_0$ is a square with sides of length $a$, and
let the cost functions $E,\phi$ be as in (\ref{EP}). Assuming that the control is
applied   along the portion of the boundary with maximum curvature, we describe here the evolution of the set $\Omega(t)$.     As usual, we
denote by $B(x,r)$ the open ball centered at $x$ with radius $r$,
and by $B(\Omega_0, \rho)$ the neighborhood of radius $\rho$
around the set $\Omega_0$.
For $t, r>0$, consider the set
\bel{Vrr}
V(\Omega_0, t, r)~\doteq~\bigcup \Big\{ B(x,r)\,;~~B(x,r)\subseteq B(\Omega_0, t)\Big\}.\eeq
This is the union of all balls of radius $r$ which are entirely contained inside $B(\Omega_0, \rho)$.

\begin{figure}[ht]
\centerline{\hbox{\includegraphics[width=4cm]{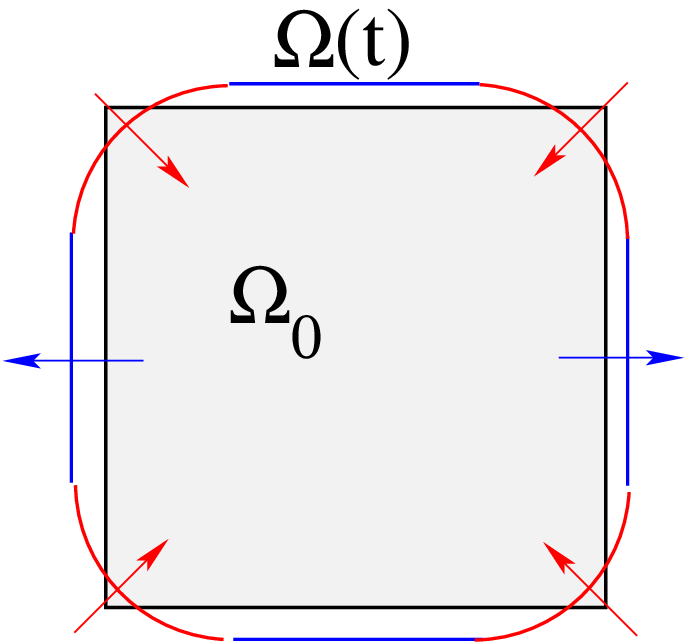}}
\qquad \hbox{\includegraphics[width=11cm]{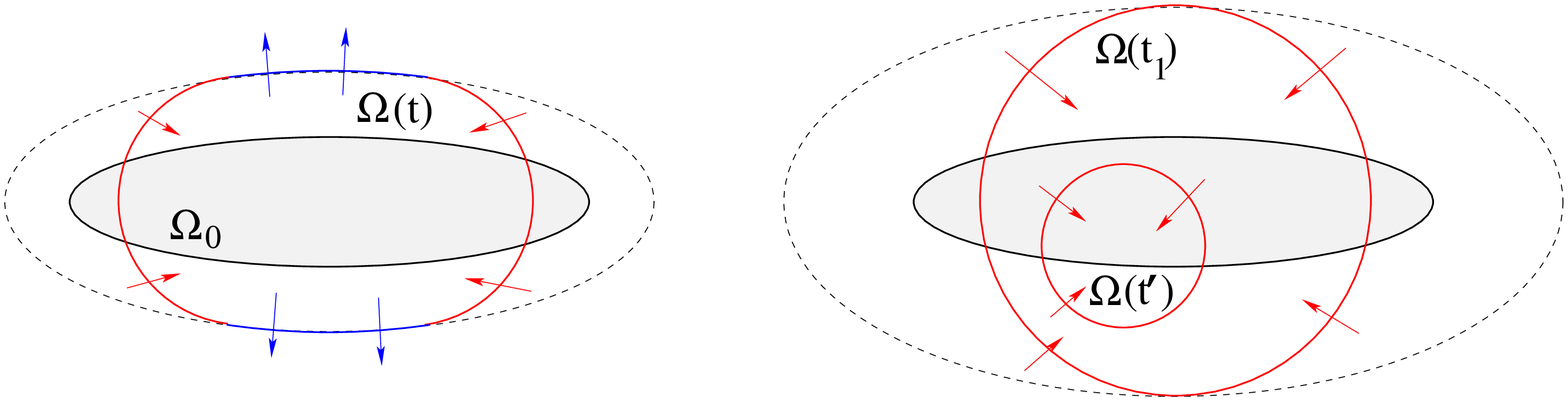}}
}
\caption{\small  Left: the moving set $\Omega(t)$ in Example~\ref{e:81}, where $\Omega_0$ is a square. Center and right:
the moving set $\Omega(t)$, in the case where $\Omega_0$ is an ellipse and the control always acts on the portion of the boundary with maximum curvature.
}
\label{f:sc62}
\end{figure}

In this case (see Fig.~\ref{f:sc62}, left), 
the set $V(\Omega_0, t,r)$ is obtained starting with a square of sides $a+2t$, 
and cutting out the regions near the four corners. The boundary of this set thus consists of four segments, and four arcs of circumferences of radius $r$.
The
perimeter and the area of the set $V(\Omega_0, t,r)$ are computed as 
\bel{PAD}\left\{ \bega{rl} P(t,r)&=~ 4(a +2t) -(8-2\pi) r ,\\[1mm]
A(t,r)&=~ (a+2t)^2- (4-\pi)r^2.\enda\right.\eeq
To derive an ODE for the function $r=r(t)$, we use the basic relation (\ref{AS}) and obtain
\bel{PAR} {d\over dt} A\bigl(\bigl(t, r(t)\bigr)-P\bigl(t, r(t)\bigr)+M~=~-2(4-\pi)r\dot{r}(t)   - (8-2\pi)r +M~=~0.
\eeq
Solving the Cauchy problem
\bel{cp1}{(8-2\pi)r\over M- (8-2\pi )r} 
\dot r ~=~ 1\,,\qquad\qquad r(0)\,=\,0,
\eeq
we determine $r(t)$  by the implicit equation
\bel{rt} \lambda \, r(t) ~=~ 1 - e^{-\lambda\bigl( t+r(t)\bigr)}\,,\qquad\qquad \lambda~\doteq~{8-2\pi\over M}.\eeq
Notice that the function $t\mapsto r(t)$ depends on $M$, but not on  $a$.

Next,  the time $t_1$ when four arcs of circumferences join together is determined by the identity
$$a+ 2t_1~=~2r(t_1).$$
This yields the implicit equation
\bel{lt1}
\lambda \left({a\over 2} + t_1\right) ~=~1 - \exp\left\{-\lambda\Big(2t_1 + {a\over 2}\Big)\right\}.\eeq
Notice that (\ref{lt1}) may not have a solution, if $M$ is too small.
To check if a solution exists set
$$f(t)~\doteq~1 - \exp\left\{-\lambda\Big(2t + {a\over 2}\Big)\right\},$$
and call 
$$\tau~=~{1\over 2\lambda} \left( \ln 2 - {a\lambda\over 2}\right)$$
the unique point where $f'(\tau) = \lambda$.
Then a solution to (\ref{lt1}) exists if and only if
\bel{tae}\lambda \left({a\over 2} + \tau\right) ~\leq ~1 - \exp\left\{-\lambda\Big(2\tau + {a\over 2}\Big)\right\}.\eeq
We observe that  a solution certainly exists if $4a<M$.  In this case, the perimeter remains strictly smaller than $M$ at all times,
and the set can be shrunk to the empty set in finite time.
\v
Finally, assuming that $t_1$ is well defined, for $t>t_1$ the optimal set $\Omega(t)$ is 
a disc whose radius satisfies
$${d\over dt}\bigl( \pi r^2(t)\bigr) ~=~ 2\pi r(t) -M.$$
\v
}\end{example}

\end{document}